# Exactness and stability in homotopical algebra (*)

MARCO GRANDIS

*Dipartimento di Matematica, Università di Genova, Via Dodecaneso 35, 16146-Genova, Italy.*
*e-mail: grandis@dima.unige.it*

**Abstract.** Exact sequences are a well known notion in homological algebra. We investigate here the more vague properties of 'homotopical exactness', appearing for instance in the fibre or cofibre sequence of a map. Such notions of exactness can be given for very general 'categories with homotopies' having *homotopy* kernels and cokernels, but become more interesting under suitable 'stability' hypotheses, satisfied – in particular – by chain complexes. It is then possible to measure the default of homotopical exactness of a sequence by the homotopy type of a certain object, a sort of 'homotopical homology'.

## Introduction

The purpose of this work is to investigate the notion of 'homotopically exact' sequence in categories equipped with homotopies, pursuing a project of developing homotopical algebra as an *enriched version* of homological algebra [G1, G2]. Well known instances of such sequences are:

(a) the *cofibre sequence*, or Puppe sequence [P1], of a map $f: A \to B$ of topological spaces or pointed spaces, where every map is, up to homotopy equivalence, an h-cokernel of the preceding one

(1) $\quad A \to B \to Cf \to \Sigma A \to \Sigma B \ ...$

($Cf$ is the h-cokernel of $f$, or standard homotopy cokernel, or mapping cone; $\Sigma$ denotes the suspension; always in the appropriate sense, absolute or pointed);

(b) the *fibre sequence* of a map $f: A \to B$ of pointed spaces, where every map is, up to homotopy equivalence, an h-kernel of the subsequent one

(2) $\quad ... \ \Omega A \to \Omega B \to Kf \to A \to B$

($Kf$ is the h-kernel, or standard homotopy kernel, or homotopy fibre; $\Omega$ denotes the loop functor);

(c) the *fibre-cofibre sequence* of a map $f: A \to B$ of chain complexes

(3) $\quad ... \ \Omega A \to \Omega B \to Kf \to A \to B \to Cf \to \Sigma A \to \Sigma B \ ...$

(*) Work supported by MURST Research Projects.

where both the aforementioned exactness conditions are satisfied, and each three-term part is homotopy equivalent to a *componentwise-split short exact sequence of complexes*.

The drastic simplification of exactness properties which appears in the last example is a product of the homotopical *stability* of chain complexes: the suspension and loop endofunctor (which in general just form an adjunction $\Sigma \dashv \Omega$), are *reciprocal* and carry the sequence (3) to itself, by a three-step shift forward $(\Sigma)$ or backward $(\Omega)$, so that the properties of its left-hand part reflect on the right-hand part and vice versa. Triangulated categories abstract these facts in the notion of 'exact triangle' (cf. Puppe [P2, P3], Verdier [Ve], Hartshorne [Ha]).

Let us go on considering this simple but relevant situation, the category of chain complexes $Ch_*\mathbf{D}$ over an *additive* category, even though the following notions are studied below in a much more general context. An *h-differential sequence* $(f, g; \alpha)$ consists of two consecutive maps of chain complexes and a nullhomotopy $\alpha$ of their composite (represented by a dotted line)

(4) $$X \xrightarrow{\ f\ } A \xrightarrow{\ g\ } Y \qquad\qquad \alpha: 0 \simeq gf: X \to Y.$$

This sequence will be said to be *h-exact* if the h-kernel of $g$ is homotopically equivalent to the h-kernel of the h-cokernel of $f$, or equivalently if the dual condition is satisfied; other conditions, of *left* and *right* h-exactness (1.3), are equivalent to the previous one in the stable case (thm. 2.4).

One can measure the default of h-exactness by the homotopy type of a complex $\mathbf{H}(f, g; \alpha)$, called the *homotopical homology* of the sequence (2.5, 3.4). It can be obtained by a *homotopical version* of the construction of ordinary homology in an abelian category $\mathbf{A}$, as presented in the left-hand diagram below for a differential sequence $(f, g)$: the construction only uses kernels and cokernels, and the sequence is exact iff $H(f, g) = 0$

(5) $$\begin{array}{ccccc} X & = & X & \longrightarrow & 0 \\ \downarrow & & \downarrow f & & \downarrow \\ Ker(g) & \rightarrowtail & A & \xrightarrow{g} & Y \\ \downarrow & & \downarrow & & \| \\ H(f, g) & \rightarrowtail & Cok(f) & \longrightarrow & Y \end{array} \qquad\qquad \begin{array}{ccccc} X & = & X & & \\ u_\alpha \downarrow & & f \downarrow & \alpha & \\ \gamma u_\alpha \quad Kg & \xrightarrow{kg} & A & \xrightarrow{g} & Y \\ \downarrow cu_\alpha & & \downarrow cf & & \| \\ \mathbf{H}(f, g; \alpha) & \longrightarrow & Cf & \xrightarrow{v_\alpha} & Y \\ & & \kappa v_\alpha & & \end{array}$$

Similarly, if $\mathbf{A} = Ch_*\mathbf{D}$, replacing (co)kernels with standard *homotopy* (co)kernels, we may construct the right-hand diagram above: we start now from the h-differential sequence $(f, g; \alpha)$, construct the h-kernel $Kg$ (with a structural nullhomotopy $\kappa g: 0 \simeq g.kg$), the h-cokernel $Cf$ (with $\gamma f: 0 \simeq cf.f$), two canonical maps $u_\alpha: X \to Kg$ and $v_\alpha: Cf \to Y$ (coherently with all previous maps *and* homotopies), and, rather surprisingly, find that there is *one* chain complex which is the h-kernel of $Cf \to Y$ and the h-cokernel of $X \to Kg$, namely

(6) $$\mathbf{H}_n(f, g; \alpha) = X_{n-1} \oplus A_n \oplus Y_{n+1}, \qquad \partial(x, a, y) = (-\partial x, -fx + \partial a, -\alpha x + ga - \partial y).$$

Now, the sequence $(f, g; \alpha)$ is h-exact iff $\mathbf{H}(f, g; \alpha)$ is contractible. Note that, *in contrast with ordinary homology*, all this requires the *additive* structure of $\mathbf{D}$, instead of kernels and cokernels (i.e., exactness in the sense of Puppe-Mitchell); moreover, monics are replaced with *fibrations* and epis with *cofibrations*.

The abstract frame we will use for 'categories with homotopies' is a notion of *homotopical category*, developed in previous papers and recalled in 1.1-2: it is a sort of lax 2-category with suitable comma and cocomma squares. The main reference, [G2], is cited here as Part I (and I.7 means Section 7 therein). The links of this setting with triangulated categories are dealt with in [G4]; its 'parallelism' with homological algebra in [G1]. Finally, it should be noted that the Freyd embedding of a stable homotopy category into an abelian category [F1, F2, F3] can be extended to the unstable case [G6], providing a notion of exactness for a sequence (f, g) of maps with $gf \simeq 0$, which is *weaker* than h-exactness and only depends on the homotopy classes [f], [g] (1.3, 2.2).

*Outline*. Section 1 introduces our notions of homotopical exactness of an h-differential sequence. Then, in Section 2, various properties of stability are considered (2.2-3), with their links with h-exactness of sequences; e.g., in the h-stable case, all our notions of homotopical exactness become equivalent (2.4) and their default is measured by the homotopy type of suitable objects (2.5).

Section 3 considers h-exactness for the categories of chain complexes (a stable case), positive complexes (left h-stable) and bounded complexes between fixed degrees, $\mathrm{Ch}_0^p\mathbf{D}$ (h-semistable); all these homotopical categories have a homotopical homology, which characterises h-exactness in the unbounded and positive cases (3.4, 3.6). In the bounded case $\mathrm{Ch}_0^p\mathbf{D}$, we have similar results up to weak equivalences (the chain maps which induce isomorphism in homology): the *weak* homotopy type of $\mathbf{H}(f, g; \alpha)$ measures the *weak* exactness of the sequence (3.7-9). Finally, the case $p = 1$ is considered in Section 4: $\mathrm{Ch}_0^1\mathbf{D} = \mathbf{D^2}$ is the category of maps of $\mathbf{D}$, with suitable homotopies; here, homotopical homology does measure the default of h-exactness, which is equivalent to weak exactness. In particular, as a special case of the Dold-Kan theorem [Ka], $\mathbf{Ab^2}$ is equivalent to $\mathbf{Cat}(\mathbf{Ab})$, the category of internal categories in abelian groups, and the present h-exactness property agrees with a notion of exactness in *symmetric cat-groups* recently introduced by Kasangian-Vitale [KV].

## 1. Exactness in homotopical categories

After a brief review of homotopical categories (1.1-2) from Part I [G2], we define some notions of homotopical exactness for homotopically differential sequences (1.3), and prove that the fibre and cofibre sequences of a map are, respectively, left and right h-exact (1.7-8).

**1.1. Homotopical categories.** $\mathbf{A}$ will always be a *pointed homotopical category*, as defined in I.7. The precise definition is too long to be rewritten here, we just review its main points, also to fix the (slightly different) notation we use here.

To begin with, $\mathbf{A}$ is a sort of lax 2-category, quite different from bicategories. It has objects, maps, and *homotopies* (2-cells) $\alpha: f \simeq g: A \to B$, with much of the usual structure of 2-dimensional categories, plus an equivalence relation $\alpha \simeq_2 \alpha'$ between parallel homotopies $\alpha, \alpha': f \simeq g$, called 2-homotopy. The vertical composition (or *concatenation*) of homotopies will be written additively, $\alpha + \beta: f \simeq h$ (for $\beta: g \simeq h$); a vertical identity (or *trivial homotopy*) as $0: f \simeq f$ or $0_f$; a *reverse* homotopy as $-\alpha: g \simeq f$. But, of course, the vertical structure *only* behaves categorically ('groupoidally') *up to 2-homotopy*. There is only a *reduced* horizontal composition of homotopies and maps $v.\alpha.u$ (*whisker composition*, for $u: A' \to A$, $v: B \to B'$); the *reduced interchange property* holds up to 2-homotopy (which is why there is no general horizontal composition of 2-cells).

The homotopy relation $f \simeq g$ (meaning that there is some homotopy $\alpha: f \simeq g$) is a congruence of categories; the quotient $\mathrm{Ho}\mathbf{A} = \mathbf{A}/\!\simeq$ is called the *homotopy category* of $\mathbf{A}$ (and coincides, under mild supplementary hypotheses, with the category of fractions which inverts homotopy equivalences, in the obvious sense). The 2-category $\mathrm{Ho}_2\mathbf{A} = \mathbf{A}/\!\simeq_2$ (with the same objects, same maps, and *tracks*, i.e. classes of homotopies up to 2-homotopy) is called the *track 2-category* of $\mathbf{A}$.

It is important to note that, generally, one cannot reduce the study of homotopy to $\mathrm{Ho}_2\mathbf{A}$, somewhat in the same way as higher dimensional category theory cannot be reduced to 2-categories. Homotopies can usually be represented by the cylinder functor, or dually by the path functor; often by both, via their adjunction $I \dashv P$; then, all *higher homotopies* are automatically produced by their powers. This gives a more powerful abstract frame, studied in [G5] in a form – IP4-homotopical categories – which will be marginally used here (even if the comparison of the left and right homotopical homologies, in 2.5, requires such a stronger setting).

**1.2. Homotopy kernels.** As a second main aspect, the pointed homotopical category $\mathbf{A}$ is assumed to have homotopy kernels and homotopy cokernels, with respect to a zero object 0; the latter is defined by a 2-dimensional universal property: every object $X$ has precisely one map $t: X \to 0$ and one homotopy $t \simeq t$ (necessarily $0_t$); and dually.

Also to fix the present notation, the *standard homotopy kernel*, or *h-kernel*, of the map $f: A \to B$ is a triple $\mathrm{hker}(f) = (Kf, kf, \kappa f)$, as in the left-hand diagram below, determined up to isomorphism by the following universal property (of comma squares)

(1)
$$\begin{array}{ccc} A & \xrightarrow{f} & B \\ {\scriptstyle kf}\uparrow & \nwarrow {\scriptstyle \kappa f} & \uparrow \\ Kf & \longrightarrow & 0 \end{array} \qquad\qquad \begin{array}{ccc} 0 & \longrightarrow & B \\ \uparrow & \nwarrow {\scriptstyle \omega B} & \uparrow \\ \Omega B & \longrightarrow & 0 \end{array}$$

- for every similar triple $(X, x, \xi)$, where $x: X \to A$ and $\xi: 0 \simeq fx: X \to B$, there is a unique $u: X \to Kf$ such that $kf.u = x$ and $\kappa f.u = \xi$.

In particular, for $A = 0$ (as in the right-hand diagram above), we obtain the *loop-object* $\Omega B = K(0 \to B)$, with a structural homotopy $\omega_B: 0 \simeq 0: \Omega B \to B$.

The dual universal property determines the *h-cokernel* $\mathrm{hcok}(f) = (Cf, cf, \gamma f)$

(2)
$$\begin{array}{ccc} A & \xrightarrow{f} & B \\ \downarrow & {\scriptstyle \gamma f} \nearrow & \downarrow {\scriptstyle cf} \\ 0 & \longrightarrow & Cf \end{array} \qquad\qquad \begin{array}{ccc} A & \longrightarrow & 0 \\ \downarrow & {\scriptstyle \sigma A} \nearrow & \downarrow \\ 0 & \longrightarrow & \Sigma A \end{array}$$

which reduces to the *suspension* $\Sigma A$, when $B = 0$. *Both constructs are also assumed to satisfy a 2-dimensional universal property* (I.2.6), which we do not recall here but will be used in some proofs (e.g. in 1.6; 2.5) and through various results of Part I which depend on it. One of its consequences is the fact that each h-kernel map $kf$ is a *fibration* (I.3.6), and each h-cokernel map $cf$ is a *cofibration*.

The h-kernel of $A \to 0$ is $1: A \to A$ (with the trivial endohomotopy of $A \to 0$) and the h-cokernel of $0 \to A$ is $1: A \to A$. But note that $f \simeq 0$ (or even $f = 0$) does not imply that $kf: Kf \to A$ (or $cf$) be a homotopy equivalence, as shown by $\Omega B = K(0 \to B)$.

An object $X$ is said to be *contractible*, or *nullhomotopic*, if it is homotopy equivalent to 0, or equivalently if $1_X \simeq 0\colon X \to X$. If $f\colon A \to B$ is a homotopy equivalence, then $Kf$ and $Cf$ are contractible, since the reflection theorem for h-pullbacks (I.3.7) shows that $Kf \to 0$ is also a homotopy equivalence (the converse holds under some stability hypotheses, 2.4).

Thus, for every object $A$, the *cocone* $KA = K(1_A)$ and the *cone* $CA = C(1_A)$ are contractible

$$
(3) \qquad
\begin{array}{ccc}
A & = & A \\
{\scriptstyle kA}\uparrow & \nwarrow {\scriptstyle \kappa A} & \uparrow \\
KA & \longrightarrow & 0
\end{array}
\qquad\qquad
\begin{array}{ccc}
A & = & A \\
\downarrow & {\scriptstyle \gamma A} \nearrow & \downarrow {\scriptstyle cA} \\
0 & \longrightarrow & CA
\end{array}
$$

**1.3. Homotopical exactness.** An *h-differential* sequence $(f, g; \alpha)$ consists of a nullhomotopy $\alpha\colon 0 \simeq gf$. The exactness properties we want to investigate deal with the following diagram, constructed by means of h-kernels and h-cokernels

$$
(1) \qquad
\begin{array}{ccccccc}
X & & & & & & Y \\
{\scriptstyle u_\alpha}\Big\downarrow \downarrow {\scriptstyle u_f} & \overset{f}{\searrow} & & \overset{\alpha}{- - -} & & \overset{g}{\nearrow} & {\scriptstyle v_g}\uparrow \Big\uparrow {\scriptstyle v_\alpha} \\
\bullet & \xrightarrow[\;kcf\;]{} & & A & & \xrightarrow[\;ckg\;]{} & \bullet \\
\downarrow {\scriptstyle u} & \underset{kg}{\nearrow} & & & & \underset{cf}{\searrow} & {\scriptstyle v}\uparrow \\
Kg & & & & & & Cf
\end{array}
\qquad
\begin{array}{l}
\alpha\colon\ 0 \simeq gf\colon\ X \to Y \\
u \ = \ k(v_\alpha), \quad v \ = \ c(u_\alpha),
\end{array}
$$

the vertical arrows are defined via the universal properties, by the following conditions (the assertion $u = k(v_\alpha)$, $v = c(u_\alpha)$ is motivated by the functors $c\colon \mathbf{A}/A \rightleftarrows \mathbf{A}\backslash A :k$ described below, 1.4.5)

$$
\begin{array}{llllll}
(2) & kg.u_\alpha \ = \ f, & \kappa g.u_\alpha \ = \ \alpha, & & v_\alpha.cf \ = \ f, & v_\alpha.\gamma f \ = \ \alpha, \\
& kcf.u_f \ = \ f, & \kappa cf.u_f \ = \ \gamma f, & & v_g.ckg \ = \ g, & v_g.\gamma kg \ = \ \kappa g, \\
& kg.u \ = \ kcf, & \kappa g.u \ = \ v_\alpha.\kappa cf, & & v.cf \ = \ ckg, & v.\gamma f \ = \ \gamma kg.u_\alpha.
\end{array}
$$

(Note that $u_f = u_{\gamma f}$, $v_f = v_{\kappa f}$). The diagram is commutative, i.e. $u.u_f = u_\alpha$ and $v_g.v = v_\alpha$, since

$$
(3) \quad kg.(uu_f) \ = \ kcf.u_f \ = \ f, \qquad\qquad \kappa g.(uu_f) \ = \ v_\alpha.\kappa cf.u_f \ = \ v_\alpha.\gamma f \ = \ \alpha.
$$

We say that the h-differential sequence $(f, g; \alpha)$ is:

(a) *left h-exact* if the map $u_\alpha\colon X \to Kg$, determined by the universal property of the h-kernel, is a homotopy equivalence;

(a*) *right h-exact* if $v_\alpha\colon Cf \to Y$ is a homotopy equivalence;

(b) *strongly h-exact* if it is both left and right h-exact;

(c) *h-exact* if $u = k(v_\alpha)\colon Kcf \to Kg$ and $v = c(u_\alpha)\colon Cf \to Ckg$ are homotopy equivalences.

By (2), $u = u_\lambda$ and $v = v_\mu$, with respect to the nullhomotopies $\lambda = v_\alpha.\kappa cf$ and $\mu = \gamma kg.u_\alpha$

$$
(4) \qquad Kcf \xrightarrow{\;kcf\;} A \xrightarrow{\;g\;} Y \ \ (\lambda) \qquad\qquad X \xrightarrow{\;f\;} A \xrightarrow{\;ckg\;} Ckg \ \ (\mu)
$$

whence $(f, g; \alpha)$ is h-exact iff the sequences above are *left* and *right* h-exact, respectively.

We shall see that strong h-exactness implies h-exactness (end of 1.5). Strong h-exactness is often too strong to be of interest (e.g. for pointed spaces), but is of use for stable homotopical categories, where all these conditions are equivalent, or also under weaker stability conditions (2.4); the last condition becomes simpler in the 'semistable' case, where its two requirements are equivalent (2.4d) and therefore implied both by left or right h-exactness.

It is interesting to note what happens when the middle object A is zero

$$
(5) \qquad
\begin{array}{ccccccc}
X & & & & & & Y \\
\downarrow u_X & \searrow & & \overset{\alpha}{- - - -} & & \nearrow & \uparrow v_Y \\
u_\alpha \; \Omega\Sigma X & \longrightarrow & & 0 & & \longrightarrow & \Sigma\Omega Y \; v_\alpha \\
\downarrow u & \nearrow & & & & \searrow & \uparrow v \\
\Omega Y & & & & & & \Sigma X
\end{array}
\qquad
\begin{array}{l}
\alpha\colon 0 \simeq 0\colon X \to Y \\
u = \Omega(v_\alpha), \quad v = \Sigma(u_\alpha)
\end{array}
$$

in this case (which would be trivial in homological algebra), our data reduce to a triple (X, Y; α), formed of two objects and a nullhomotopy $\alpha\colon 0 \simeq 0\colon X \to Y$. The typical left h-exact sequence of this type is $(\Omega Y, Y; \omega_Y)$, and every other is coherently equivalent to one of this type, in the sense of definition 1.6; similarly, the typical right h-exact sequence is $(X, \Sigma X; \sigma_X)$. This case also shows that a right h-exact sequence $(X, \Sigma X; \sigma_X)$ need not be h-exact: in the homotopical category **Top•** of pointed spaces, take $X = \mathbf{S}^0$; then $\Sigma X = \mathbf{S}^1$ is not homotopy equivalent to $\Sigma\Omega\mathbf{S}^1$ ($H_1(\Sigma\Omega\mathbf{S}^1) = \tilde{H}_0(\Omega\mathbf{S}^1)$ is a free abelian group of countable rank).

Marginally, we also consider a condition of *pseudo-exactness*, meaning that $gf \simeq 0$ and $cf.kg \simeq 0$, which actually just concerns a pair (f, g) of consecutive arrows instead of a sequence (f, g; α); any h-exact sequence is also pseudo-exact, since in (1) we have $cf.kg.u = cf.kcf \simeq 0$ (and u is assumed to be a homotopy equivalence), but the converse is far from being true: any sequence $X \to 0 \to Y$ is trivially pseudo-exact, but – for any choice of α – *cannot be left h-exact unless* X is homotopy equivalent to ΩY (and h-exact coincides with left h-exact, in the stable case, by 2.4).

**1.4. Slice categories.** The kernel-cokernel adjunction of abelian categories, between slice categories, corresponds here to an adjunction given by h-kernels and h-cokernels, which we have already used, implicitly, to define h-exactness. We consider first its strict version, then – in the next subsection – its coherent version.

By h-kernels and h-cokernels, a morphism $f = (f', f'')\colon x \to y$ of the category $(\mathbf{2}, \mathbf{A}) = \mathbf{A}^\mathbf{2}$ of maps of **A** yields a commutative diagram, where K(f) and C(f) are defined as follows

$$
(1) \qquad
\begin{array}{ccccccc}
Kx & \xrightarrow{kx} & X' & \xrightarrow{x} & X'' & \xrightarrow{cx} & Cx \\
{\scriptstyle K(f)}\downarrow & & {\scriptstyle f'}\downarrow & & \downarrow{\scriptstyle f''} & & \downarrow{\scriptstyle C(f)} \\
Ky & \xrightarrow[ky]{} & Y' & \xrightarrow[y]{} & Y'' & \xrightarrow[cy]{} & Cy
\end{array}
$$

(2) $ky.K(f) = f'.kx, \qquad \kappa y.K(f) = f''.\kappa x,$

$\phantom{(2)}$ $C(f).cx = cy.f'', \qquad C(f).\gamma x = \gamma y.f'.$

This gives two adjoint endofunctors $c \dashv k\colon \mathbf{A}^\mathbf{2} \to \mathbf{A}^\mathbf{2}$, whose action on maps is obvious

(3) $k(f) = (K(f), f')\colon kx \to ky, \qquad c(f) = (f'', C(f))\colon cx \to cy,$

and will always be written with parentheses, *to avoid ambiguity with the structural morphisms* cx, ky. The unit and counit are again provided by the universal properties

(4) $(u_x, 1): x \to kcx$, $\qquad$ $(1, v_y): cky \to y$,

$$\begin{array}{ccccc} \bullet & \xrightarrow{x} & A & \xrightarrow{cx} & Cx \\ {\scriptstyle u_x}\downarrow & \nearrow{\scriptstyle kcx} & & & \\ Kcx & & & & \end{array} \qquad \begin{array}{ccccc} Kx & \xrightarrow{ky} & A & \xrightarrow{y} & \bullet \\ & & & {\scriptstyle cky}\searrow & \uparrow{\scriptstyle v_y} \\ & & & & Cky \end{array}$$

$kcx.u_x = x, \quad \kappa cx.u_x = \gamma(x),$ $\qquad$ $v_y.cky = y, \quad v_y.\gamma ky = \kappa y.$

Keeping fixed the object A above, all this restricts to an adjunction between the slice categories **A**/A = (**A**↓A) of *objects over* A, x: • → A, and **A**\A = (A↓**A**) of *objects under* A, y: A → •

(5) $c: \mathbf{A}/A \rightleftarrows \mathbf{A}\backslash A : k$, $\qquad$ $u_x: x \to kcx, \quad v_y: cky \to y$,

$$\begin{array}{ccccc} X & \xrightarrow{x} & A & \xrightarrow{cx} & Cx \\ {\scriptstyle f}\downarrow & & \| & & \downarrow{\scriptstyle c(f)} \\ X' & \xrightarrow[x']{} & A & \xrightarrow[cx']{} & Cx' \end{array} \qquad \begin{array}{ccccc} Kx & \xrightarrow{ky} & A & \xrightarrow{y} & Y \\ {\scriptstyle k(f)}\downarrow & & \| & & \downarrow{\scriptstyle f} \\ Ky & \xrightarrow[ky']{} & A & \xrightarrow[y]{} & Y' \end{array}$$

For A = 0, both slice categories reduce to **A** and the previous adjunction (5) becomes the classical suspension-loop adjunction

(6) $\Sigma: \mathbf{A} \rightleftarrows \mathbf{A} : \Omega$, $\qquad$ $u_X: X \to \Omega\Sigma(X), \quad v_X: \Sigma\Omega(X) \to X.$

All this just needs '1-dimensional homotopical properties'; precisely, it already holds in a pointed *semi*homotopical category **A** (I.7.1).

**1.5. Coherent slice categories.** But the fact that **A** is homotopical allows us to replace the strict slice category **A**/A with the more interesting *coherent homotopy category over* A (cf. [HK1-2], [HKP]), which will be written [**A**/A].

To begin with, we replace $\mathbf{A}^2 = (\mathbf{2}, \mathbf{A})$ with the *coherent homotopy category of morphisms* $[\mathbf{2}, \mathbf{A}]$: an object is still an **A**-map $x: X' \to X''$, but a morphism $[f] = [f', f''; \varphi]: x \to y$ is determined by a triple $f = (f', f''; \varphi)$ forming a homotopy-commutative square $\varphi: f''x \simeq yf'$ as in the left-hand diagram below

(1)
$$\begin{array}{ccc} X' & \xrightarrow{x} & X'' \\ {\scriptstyle f'}\downarrow & \varphi \swarrow & \downarrow{\scriptstyle f''} \\ Y' & \xrightarrow[y]{} & Y'' \end{array} \qquad \begin{array}{ccccc} X' & = & X' & \xrightarrow{x} & X'' \\ {\scriptstyle g'}\downarrow & \xleftarrow{\alpha'} & \downarrow{\scriptstyle f'} & \varphi\swarrow & \downarrow{\scriptstyle f''} \\ Y' & = & Y' & \xrightarrow[y]{} & Y'' \end{array} \quad \simeq_2 \quad \begin{array}{ccccc} X' & \xrightarrow{x} & X'' & = & X'' \\ {\scriptstyle g'}\downarrow & \swarrow\psi & {\scriptstyle g''}\downarrow & \xleftarrow{\alpha''} & \downarrow{\scriptstyle f''} \\ Y' & \xrightarrow[y]{} & Y'' & = & Y'' \end{array}$$

up to identifying $[f', f''; \varphi] = [g', g''; \psi]$ when there exists a *coherent* pair of homotopies $\alpha': f' \to g'$, $\alpha'': f'' \to g''$, in the sense that $\varphi + y\alpha' \simeq_2 \alpha''x + \psi$ (as in the right-hand diagram above). The composition is obtained from the pasting of homotopies, and yields indeed a category (since we are actually working in a 2-category with invertible cells, $Ho_2\mathbf{A}$).

As a crucial fact, well known for the coherent category of objects over a fixed space ([HK2], thm. 1.3), and related to a classical theorem of Dold ([Do], 6.1), the map $[f', f''; \varphi]: x \to y$ is an isomorphism iff $f'$ and $f''$ are homotopy equivalences. In fact, in this case, one can choose an *h-adjoint equivalence* for $f'$

(2) $\alpha': 1 \to g'f', \quad \beta': f'g' \to 1 \qquad\qquad (f'\alpha' + \beta'f' \simeq_2 0_{f'}, \quad \alpha'g' + g'\beta' \simeq_2 0_{g'}),$

(by Vogt's Lemma [Vo]: given an arbitrary equivalence, replace $\beta'$ with $(-\beta'f'g' - f'\alpha'g') + \beta'$ and verify the triangle identities); and similarly for $f''$. Finally, one constructs an inverse $[g', g'', \psi]$ with a suitable homotopy $\psi: g''y \to xg'$ (namely, $\psi = (-g''y\beta' - g''\varphi g') - \alpha''xg'$)

Replacing diagram 1.4.1 with the following one (where $f = (f', f''; \varphi)$)

(3)
$$\begin{array}{ccccccc}
Kx & \xrightarrow{kx} & X' & \xrightarrow{x} & X'' & \xrightarrow{cx} & Cx \\
{\scriptstyle K(f)}\downarrow & & {\scriptstyle f'}\downarrow & \overset{\varphi}{\swarrow} & \downarrow{\scriptstyle f''} & & \downarrow{\scriptstyle C(f)} \\
Ky & \xrightarrow[ky]{} & Y' & \xrightarrow[y]{} & Y'' & \xrightarrow[cy]{} & Cy
\end{array}$$

(4) $ky.K(f) = f'.kx, \qquad\qquad \kappa y.K(f) = f''.\kappa x + \varphi.kx: \ 0 \simeq f''x.kx \to yf'.kx,$

$C(f).cx = cy.f'', \qquad\qquad C(f).\gamma x = \gamma y.f' - cy.\varphi: 0 \simeq cy.yf' \to cy.f''x,$

we have now two endofunctors $c, k: [\mathbf{2}, \mathbf{A}] \to [\mathbf{2}, \mathbf{A}]$

(5) $k[f] = [K(f), f'; 0]: kx \to ky, \qquad\qquad c[f] = [f'', C(f); 0]: cx \to cy,$

which can be easily shown to be indeed well defined. They are adjoint, $c \dashv k$, with (strict) unit and counit defined essentially as above ($u_x$ and $v_y$ are defined by the same formulae, 1.4.4)

(6) $[u_x, 1; 0]: x \to kcx, \qquad\qquad [1, v_y; 0]: cky \to y.$

For a fixed object A, we have a restricted adjunction between *h-objects over* A, $x: \bullet \to A$, and *h-objects under* A, $y: A \to \bullet$

(7) $c: [\mathbf{A}/A] \rightleftarrows [\mathbf{A}\backslash A] : k, \qquad\qquad [u_x]: x \to kcx, \quad [v_y]: cky \to y;$

Now, the category $[\mathbf{A}/A] \subset [\mathbf{2}, \mathbf{A}]$ has morphisms $[f, \varphi] = [f, 1_A; \varphi]: x \to y$ determined by a homotopy-commutative triangle $\varphi: yf \simeq x$; similarly $[\mathbf{A}\backslash A] \subset [\mathbf{2}, \mathbf{A}]$. For $A = 0$, this reduces again to a suspension-loop adjunction, now *at the level of the homotopy category* $\mathrm{Ho}\mathbf{A} = [\mathbf{A}/0] = [\mathbf{A}\backslash 0]$, and induced by the previous one (1.4.6)

(8) $\Sigma: \mathrm{Ho}\mathbf{A} \rightleftarrows \mathrm{Ho}\mathbf{A} : \Omega, \qquad\qquad [u_X]: X \to \Omega\Sigma(X), \quad [v_X]: \Sigma\Omega(X) \to X.$

Thus, in 1.3, the h-differential sequence $(f, g; \alpha)$ is *left h-exact* iff the map $[u_\alpha]: f \to kg$ is an isomorphism of $[\mathbf{A}/A]$; and *h-exact* iff $k[v_\alpha]: kcf \to kg$ *and* $c[u_\alpha]: cf \to ckg$ are isomorphisms of $[\mathbf{A}/A]$ and $[\mathbf{A}\backslash A]$, respectively. This makes evident that strong h-exactness implies h-exactness.

**1.6. Equivalence of sequences.** Given a coherent diagram, where the vertical arrows are homotopy equivalences

(1)
$$\begin{array}{ccccc} & f & \overset{\alpha}{- - - -} & g & \\ A' & \longrightarrow & A & \longrightarrow & A'' \\ u' \downarrow & \varphi \swarrow & \downarrow u & \psi \swarrow & \downarrow u'' \\ B' & \longrightarrow & B & \longrightarrow & B'' \\ & x & \underset{\beta}{- - - -} & y & \end{array} \qquad \begin{array}{l} \alpha: 0 \simeq gf, \qquad \beta: 0 \simeq yx \\ \varphi: uf \simeq xu', \qquad \psi: u''g \simeq yu \\ \beta u' - y\varphi \simeq_2 u''\alpha + \psi f, \end{array}$$

we say that the sequences $(f, g; \alpha)$ and $(x, y; \beta)$ are (coherently) *equivalent*. Plainly, a left h-exact sequence $(f, g; \alpha)$ is equivalent to $(kg, g; \kappa g)$. Also the converse holds: more generally, *all the h-exactness conditions considered in* 1.3 *are invariant up to equivalence*.

In fact, let us assume that the upper row is left h-exact and prove that the lower one is too. (The similar fact for h-exactness follows then from the characterisation of the latter in 1.3.4). The map $[u, u''; \psi]: g \to y$ is iso in $[\mathbf{2}, \mathbf{A}]$, whence $k[u, u''; \psi]$ is iso, and $w = K(u, u''; \psi)$ is a homotopy equivalence. We construct the diagram

(2)
$$\begin{array}{ccccccc} & s & & kg & \overset{\kappa g}{- - - -} & g & \\ A' & \longrightarrow & Kg & \longrightarrow & A & \longrightarrow & A'' \\ u' \downarrow \ \swarrow \varphi' & & w \downarrow & & \downarrow u & \psi \swarrow & \downarrow u'' \\ B' & \longrightarrow & Ky & \longrightarrow & B & \longrightarrow & B'' \\ & t & & ky & \underset{\kappa y}{- - - -} & y & \end{array}$$

(3) $ky.w = u.kg, \qquad \kappa y.w = u''.\kappa g + \psi.kg,$

$kg.s = f, \qquad \kappa g.s = \alpha, \qquad\qquad ky.t = x, \qquad \kappa y.t = \beta,$

where $s = u_\alpha$ is a homotopy equivalence, by hypothesis, and we want to prove that also $t = u_\beta$ is so; plainly, it suffices to prove the existence of a homotopy $\varphi': ws \to tu'$, which can be derived from the 2-dimensional universal property of the h-pullback $hker(y) = (Ky, ky, \kappa y)$. In fact, the morphisms $ws, tu': A' \to Ky$ have a homotopic projection on $B$, coherently with $\kappa y: 0 \simeq y.ky$ (as a consequence of the coherence hypothesis, in (1))

(4) $ky.ws = u.kg.s = uf, \qquad ky.tu' = xu', \qquad \varphi: ky.ws \to ky.tu',$

$\kappa y.ws + y\varphi = (u''.\kappa g + \psi.kg).s + y\varphi = (u''\alpha + \psi f) + y\varphi \simeq_2 \beta u' = \kappa y.tu'.$

**1.7. Theorem: exactness of the cofibre sequence.** If $\mathbf{A}$ is pointed homotopical, the cofibre sequence of a map $f: A \to B$

(1)
$$A \xrightarrow{f} B \xrightarrow{cf} Cf \xrightarrow{\delta} \Sigma A \xrightarrow{\Sigma f} \Sigma B \xrightarrow{\Sigma cf} \Sigma Cf \ \ldots$$
$$\qquad \underset{\gamma f}{- - - -} \qquad \underset{0}{- - - -} \qquad \underset{\sigma f}{- - - -} \qquad \underset{\gamma' f}{- - - -}$$

is right h-exact at any point, with non-canonical nullhomotopies $\sigma f: 0 \simeq \Sigma f.\delta$, $\gamma' f: 0 \simeq \Sigma cf.\Sigma f$, and so on. In the stronger assumption of a pointed IP4-homotopical structure for $\mathbf{A}$ [G5], there are canonical nullhomotopies $\sigma f$, $\Sigma(\gamma f)$,... having that result.

**Proof.** The first part of the statement is essentially proved in I.5 (and only needs a *right*-homotopical category). To begin with, it is proved (I.5.7.11) that there exists *some* nullhomotopy $\sigma f$ satisfying

(2) $\sigma f: 0 \simeq \Sigma f.\delta: Cf \to \Sigma B, \qquad\qquad \sigma f.cf = \sigma_B: 0 \simeq 0: B \to B.$

Then, applying the homotopy invariant functor $\Sigma$ (I.4.5) we obtain, again in a non-canonical way, the subsequent nullhomotopies $\gamma'f$,... We have now the *contracted cofibre diagram* (I.5.6.3)

$$(3)\qquad \begin{array}{ccccccccccccc} A & \xrightarrow{f} & B & \xrightarrow{cf} & Cf & \xrightarrow{\delta} & \Sigma A & \xrightarrow{\Sigma f} & \Sigma B & \xrightarrow{\Sigma cf} & \Sigma Cf & \ldots \\ \| & & \| & & \| & & \uparrow u_0 & \simeq & \uparrow u_1 & \simeq & \uparrow u_2 & \\ A & \xrightarrow[f]{} & B & \xrightarrow[x_1]{} & B_2 & \xrightarrow[x_2]{} & B_3 & \xrightarrow[x_3]{} & B_4 & \xrightarrow[x_4]{} & B_5 & \ldots \end{array}$$

linking the cofibre sequence of $f$ to the sequence of its iterated h-cokernels

(4) $x_1 = cf,\qquad x_{i+1} = cx_i,\quad \mathrm{hcok}(x_i) = (x_{i+1}, \gamma x_i\colon 0 \simeq x_{i+1}.x_i)\qquad (i \geq 1);$

all squares are h-commutative, all vertical arrows are homotopy equivalences. Moreover, the diagram is coherent (as in 1.6.1), with respect to the nullhomotopies of the upper row (as in (1)) and the structural nullhomotopies $\gamma x_i$ of the lower one: this follows from the construction of $u_i$'s in I.5.4-6.

Now, if **A** is pointed I4-homotopical, the cone functor $C\colon \mathbf{A} \to \mathbf{A}$ inherits an induced monad (as proved more in detail in [G3, 3.7]),

(5) $c\colon 1 \to C,\qquad \mathbf{g}\colon C^2 \to C \qquad (\mathbf{g}.cC = \mathrm{id} = \mathbf{g}.Cc),$

whose operation **g** (induced by a 'connection' of the cylinder functor) represents a *natural* nullhomotopy of the cone $CX$, $\mathbf{g}X.\gamma CX\colon 0 \simeq 1\colon CX \to CX$. We can now deduce a structural nullhomotopy $\sigma f\colon 0 \simeq \Sigma f.\delta$; in fact, $\Sigma f.\delta$ factors through $CB$, as proved by the following (canonical) diagram

(6)
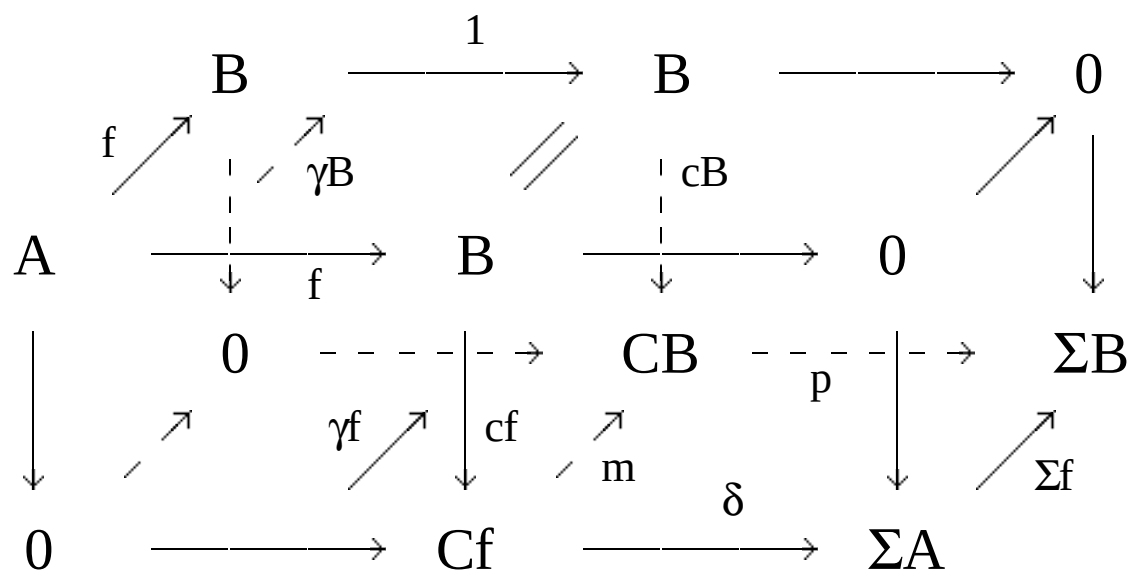

In the left-hand cube, the front and back face are h-pushouts (the h-cokernels of $f$ and $1_B$); in the right-hand cube, the front and back face are ordinary pushouts, so that the front and back rectangle are h-pushouts, by the pasting property I.2.2: $\sigma_A = \delta.\gamma f$ and $\sigma_B = p.\gamma B$; thus $\Sigma f.\delta = (Cf \to CB \to \Sigma B) = pm$, and we can define

(7) $\sigma f = p.\mathbf{g}B.\gamma CB.m\colon 0 \simeq pm,$

$\sigma f.cf = p.\mathbf{g}B.\gamma CB.m.cf = p.\mathbf{g}B.\gamma CB.cB = p.\mathbf{g}B.CcB.\gamma B = p.\gamma B = \sigma_B.$

Finally, $\Sigma$ can be extended to homotopies (as proved for the cylinder functor, in [G5], 2.9).

**1.8. The cofibre sequence.** Dually, the fibre sequence of a map $f\colon A \to B$

$$(1)\qquad \ldots\ \Omega Kf \xrightarrow{\Omega kf} \underset{\kappa' f}{\Omega A} \xrightarrow{\Omega f} \underset{\omega f}{\Omega B} \xrightarrow{\partial} \underset{0}{Kf} \xrightarrow{kf} \underset{\kappa f}{A} \xrightarrow{f} B$$

is left-exact at any point, with non-canonical nullhomotopies $\omega f: 0 \simeq \partial.\Omega f$, $\kappa' f: 0 \simeq \Omega f.\Omega kf$, and so on. If **A** has a pointed IP4-homotopical structure [G5], these nullhomotopies can be constructed in a natural way; in particular, $\omega f$ is derived from the connection **g**: $K \to K^2$ of the cocone (the comultiplication of the comonad structure), and $\kappa' f = \Omega(\kappa f)$.

## 2. The links between stability and exactness

Some stability properties are recalled or introduced (2.2-3), by means of the previous adjunction between coherent slice categories (1.5). Under such hypotheses, the exactness properties behave in a simple way (2.4) and the default of h-exactness is measured by suitable objects (2.5).

**2.1. Exact adjunctions.** First, let us recall that an adjunction $F \dashv G$

(1) $F : \mathbf{X} \rightleftarrows \mathbf{Y} : G, \qquad u: 1_{\mathbf{X}} \to GF, \quad v: FG \to 1_{\mathbf{Y}},$

is said to be *exact* (Betti [Be]) if the natural transformations which appear in the triangle identities are isomorphisms (as it happens – automatically – for Galois connections). In other words, this means that the following equivalent conditions hold [Be, p. 46]:

(a) $Fu: F \to FGF$ is an isomorphism, (a') $vF: FGF \to F$ is an isomorphism,

(b) $uG: G \to GFG$ is an isomorphism, (b') $Gv: GFG \to G$ is an isomorphism,

(c) $FGv = vFG$, (d) $uGF = GFu$.

Of course, a reflective adjunction (i.e., $v$ is iso) is exact, as well as a coreflective one ($u$ is iso). If $F \dashv G$ is exact, as well as the consecutive adjunction $H \dashv K$

(2) $H : \mathbf{Y} \rightleftarrows \mathbf{Z} : K, \qquad u': 1 \to KH, \quad v': HK \to 1,$

*each of the following conditions is sufficient to make the composed adjunction* $HF \dashv GK$ *exact*

(i) $FGu'F = u'F$, (i*) $HKvK = vK$,

(ii) $Gu'FG = Gu'$, (ii*) $HvKH = Hv$.

For instance, assume (i) and consider the composed unit $\bar{u} = Gu'F.u: 1 \to GKHF$. Since $Fu$ and $Hu'$ are iso, also $HF\bar{u} = HFGu'F.HFu = Hu'F.HFu$ is so.

**2.2. Stability.** Homotopical stability essentially requires that the suspension-loop adjunction be an equivalence. As in Part I, we will use a stronger definition, better related with (co)fibre sequences.

In a general pointed homotopical category **A**, the fibre-cofibre sequence of a map $f: A \to B$ can be inserted as the central row of a homotopy-commutative *adjunction fibre-cofibre* diagram (I.7.5)

$$
(1)\qquad
\begin{array}{ccccccccccccccc}
 & & \dots\ \Sigma\Omega^2 B & \xrightarrow{\Sigma\Omega\partial} & \Sigma\Omega Kf & \xrightarrow{\Sigma\Omega kf} & \Sigma\Omega A & \xrightarrow{\Sigma\Omega f} & \Sigma\Omega B & \xrightarrow{\Sigma\partial} & \Sigma Kf & \xrightarrow{\Sigma kf} & \Sigma A & \xrightarrow{\Sigma f} & \Sigma B\ \dots \\
 & & {\scriptstyle v_{\Omega B}}\downarrow & & {\scriptstyle v_{Kf}}\downarrow & & \downarrow{\scriptstyle v_A} & & \downarrow{\scriptstyle v_B}\ \simeq & & \downarrow{\scriptstyle V_f}\ \simeq & & \| & & \| \\
\dots\ \Omega A & \longrightarrow & \Omega B & \xrightarrow{\partial} & Kf & \xrightarrow{kf} & A & \xrightarrow{f} & B & \xrightarrow[cf]{} & Cf & \xrightarrow[\delta]{} & \Sigma A & \longrightarrow & \Sigma B\ \dots \\
\| & & \| & & \simeq\ {\scriptstyle U_f}\downarrow & & \simeq\ \downarrow{\scriptstyle u_A} & & \downarrow{\scriptstyle u_B} & & \downarrow{\scriptstyle u_{Cf}} & & \downarrow{\scriptstyle u_{\Sigma A}} & & \\
\dots\ \Omega A & \xrightarrow[\Omega f]{} & \Omega B & \xrightarrow[\Omega cf]{} & \Omega Cf & \xrightarrow[\Omega\delta]{} & \Omega\Sigma A & \xrightarrow[\Omega\Sigma f]{} & \Omega\Sigma B & \xrightarrow[\Omega\Sigma c]{} & \Omega\Sigma Cf & \xrightarrow[\Omega\Sigma\delta]{} & \Omega\Sigma^2 A\ \dots & &
\end{array}
$$

whose upper and lower rows are obtained by letting the functors $\Sigma$ and $\Omega$ operate on the original sequence. The vertical maps consist of the units $u$ and counits $v$ of the adjunction $\Sigma \dashv \Omega$, together with the adjoint maps $U_f\colon Kf \to \Omega Cf$, $V_f\colon \Sigma Kf \to Cf$ provided by the nullhomotopy

(2) $\rho_f = cf.\kappa f - \gamma f.kf\colon\ 0 \simeq 0\colon\ Kf \to Cf.$

**A** is defined to be *h-stable* (I.7.8) if all the maps $U_f\colon Kf \to \Omega Cf$, $V_f\colon \Sigma Kf \to Cf$ are homotopy equivalences; then also the unit-maps $u_A\colon A \to \Omega\Sigma A$ and the counit-maps $v_A\colon \Sigma\Omega A \to A$ are so (because $u_A \simeq U_{(A \to 0)}$ and $v_A \simeq V_{(0 \to A)}$, I.7.6). Therefore, the endofunctors $\Sigma$ and $\Omega$ *shift the fibre-cofibre sequence of* f *forward and backward*, respectively, of three steps (up to homotopy equivalence). The homotopical category $Ch_*\mathbf{D}$ of (unbounded) chain complexes on any additive category **D** is easily seen to be *strictly stable*: all the vertical arrows of (1) are isomorphisms and all the squares commute (I.7.8; or here, 3.3).

In the h-stable case, the homotopy category Ho**A** has a canonical embedding in an abelian category Fr(Ho**A**), introduced by Freyd [F1, F2, F3] for the stable homotopy category of spaces; the same construction also works for a pointed homotopical category [G6], but provides then a pointed 'homological' category (a notion which extends Puppe-exact categories). The resulting notion of exactness, by ordinary (co)kernels, agrees with the present *pseudo-exactness* (end of 1.3): it is easy to see that the sequence (f, g) of **A** is pseudo-exact iff ([f], [g]) is exact in Fr(Ho**A**).

**2.3. Other stability properties.** We introduce now weaker conditions, related to the adjunction $c \dashv k$ between *h-objects* $f\colon \bullet \to A$ and $g\colon A \to \bullet$ (1.5). The pointed homotopical category **A** will be said to be:

(a) *left h-stable* if, for every object A, the adjunction $c\colon [\mathbf{A}/A] \rightleftarrows [\mathbf{A}\backslash A] :k$ is a coreflection, i.e. all components $u_f\colon Dom(f) \to Kcf$ are homotopy equivalences of **A**;

(a*) *right h-stable* if the dual condition holds, i.e. all components $v_g\colon Ckg \to Cod(g)$ are homotopy equivalences (if both conditions hold, the adjunction $c \dashv k$ is an adjoint equivalence);

(b) *h-semistable* if, for every object A, the adjunction $c\colon [\mathbf{A}/A] \rightleftarrows [\mathbf{A}\backslash A] :k$ is an *exact adjunction*, i.e. the following equivalent conditions hold (2.1)

- all components $c(u_f)\colon Cf \to Ckcf$ are homotopy equivalences of **A**,
- all components $v_{cf}\colon Ckcf \to Cf$ are homotopy equivalences,
- all components $u_{kg}\colon Kg \to Kckg$ are homotopy equivalences,
- all components $k(v_g)\colon Kckg \to Kg$ are homotopy equivalences.

Trivially, left (or right) h-stable implies h-semistable. Recall that the adjunction $\Sigma \dashv \Omega$ is a particular instance of the adjunction $c \dashv k$ (1.5.8); therefore, if **A** is left h-stable (resp. left and right h-stable, h-semistable), then, *at the level of the homotopy category* Ho**A**, $\Sigma \dashv \Omega$ is a coreflection (resp. an endoequivalence, an exact adjunction). There are parallel strict notions (left stable, etc.), derived from strict slice categories, where the relevant components – above – are required to be isomorphisms of **A**. ($Ch_*$**D**, being stable, is also left and right h-stable (by 2.4), but *not* strictly so (3.3). Its subcategory $Ch_p$**D** of *positive* complexes is just left h-stable (3.5), while *negative* complexes yield a right h-stable case.)

**2.4. Theorem: stability and exactness.** **A** is a pointed homotopical category.

(a) If **A** is h-stable, then it is also left and right h-stable.

(b) **A** is left h-stable if and only if all sequences $(f, cf; \gamma f)$ are strongly h-exact, if and only if every h-exact sequence is left h-exact, if and only if every right h-exact sequence is strongly h-exact. In this case, any cofibre sequence is strongly h-exact; moreover, a map f is a homotopy equivalence iff the object Cf is contractible, and this *implies* that Kf is also so.

(c) **A** is left and right h-stable if and only if all sequences $(f, cf; \gamma f)$ and $(kf, f; \kappa f)$ are strongly h-exact, if and only if all the conditions of h-exactness considered in 1.3 (left h-exact, right h-exact, strongly h-exact, h-exact sequence) are equivalent. In this case, any fibre-cofibre sequence is strongly h-exact; moreover, a map f is a homotopy equivalence iff Kf is contractible, iff Cf is so.

(d) If **A** is h-semistable, the two conditions defining h-exactness (1.3c) are equivalent: $k(v_\alpha)$ is a homotopy equivalence iff $c(u_\alpha)$ is so; whence - in this case - any left or right h-exact sequence is h-exact. (Analogously, if **A** is semistable, $k(v_\alpha)$ is iso iff $c(u_\alpha)$ is so.)

**Proof.** We use the notation of diagram 1.3.1. Recall the relations $u.u_f = u_\alpha$, $v_g.v = v_\alpha$; recall also that the fibre sequence of a map is left h-exact (1.8), and dually.

(a) Let **A** be h-stable. Because of this, the fibre-cofibre sequence of f is strongly h-exact, as follows from the adjunction fibre-cofibre diagram (2.2.1), the invariance of h-exactness (1.6), and the fact that $\Sigma$ and $\Omega$, being quasi-reciprocal, preserve left and right h-exactness. Thus, the sequences $(f, cf; \gamma f)$ and $(kf, f; \kappa f)$ are strongly h-exact, which implies the thesis, by (c).

(b) If **A** is left h-stable and the sequence $(f, g; \alpha)$ is h-exact, then $u_f: X \to Kcf$ and $u: Kcf \to Kg$ are homotopy equivalences, whence also $u_\alpha = u.u_f$ is so, which means that $(f, g; \alpha)$ is left h-exact. Trivially, if h-exact implies left h-exact, then right h-exact implies strongly h-exact. In this case, any sequence $(f, cf; \gamma f)$ is strongly h-exact. Finally, if this holds, the canonical map $u_f: X \to Kcf$ is a homotopy equivalence and **A** is left h-stable.

Now, let **A** be left h-stable. Then, we have already proved that any cofibre sequence is strongly h-exact. Secondly, if f is a homotopy equivalence, we know that Cf is contractible (1.2); conversely, in this case, the invariance of h-exactness (1.6) shows that kcf is a homotopy equivalence, whence also $f = kcf.u_f$ is so.

(c) Follows from (b), by duality.

(d) By hypothesis, $c(u_f)$ is a homotopy equivalence; if also u is so, it follows that $v = c(u_\alpha) = c(u).c(u_f)$ is a homotopy equivalence too. And dually.

**2.5. Homotopical homology.** Given an h-differential sequence $(f, g; \alpha)$, form the commutative, coherent diagram

(1) 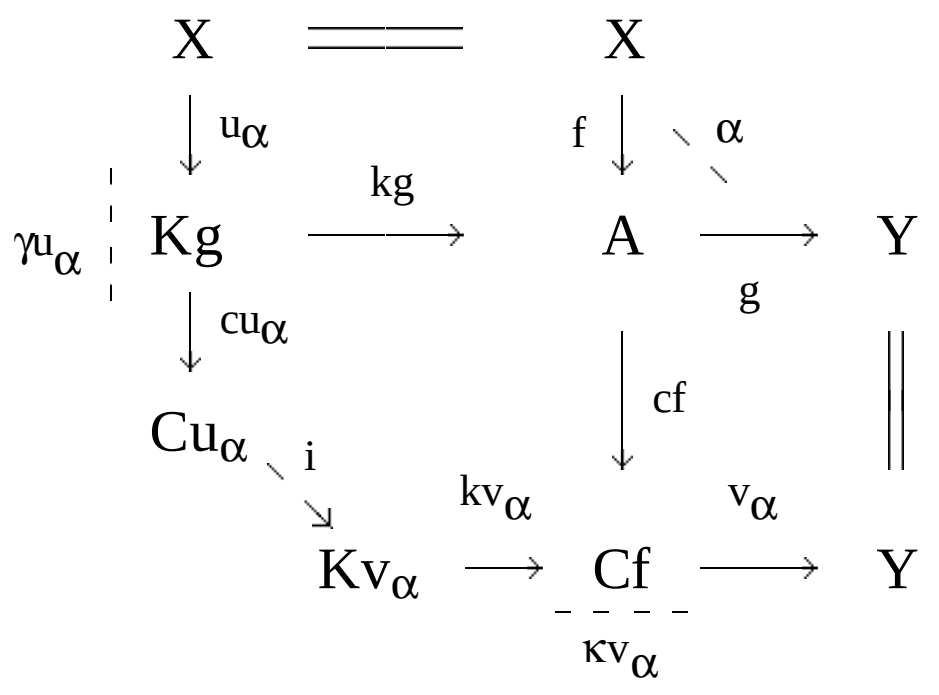

$\alpha: 0 \simeq gf: X \to Y$

$\kappa g: 0 \simeq g.kg: Kg \to Y$

$\gamma f: 0 \simeq cf.f: X \to Cf.$

It provides the *left homotopical homology* and the *right homotopical homology* of the sequence

(2) $\mathbf{H}^{-}(f, g; \alpha) = Cu_{\alpha}$, $\qquad \mathbf{H}^{+}(f, g; \alpha) = Kv_{\alpha}$,

linked by a non-canonical *comparison* $i$ satisfying the following relations (its existence is proved below; a precise determination is possible in the stronger setting of IP4-homotopical categories [G5])

(3) $kv_{\alpha}.i.cu_{\alpha} = cf.kg, \qquad kv_{\alpha}.i.\gamma u_{\alpha} = \gamma f, \qquad \kappa v_{\alpha}.i.cu_{\alpha} = \kappa g.$

In the *left h-stable case,* the homotopy type of the first object measures the default of h-exactness of $(f, g; \alpha)$: the sequence is h-exact iff it is left h-exact (2.4b), iff $u_{\alpha}$ is a homotopy equivalence, *iff* $Cu_{\alpha}$ *is contractible* (2.4b again). Dually, in the right h-stable case, $\mathbf{H}^{+}(f, g; \alpha)$ measures the same default of h-exactness, which coincides now with right h-exactness. Finally, if **A** is left and right h-stable, *both* $\mathbf{H}^{-}(f, g; \alpha)$ and $\mathbf{H}^{+}(f, g; \alpha)$ measure the default of h-exactness.

It is therefore natural to study the property, for **A**, that $\mathbf{H}^{-}(f, g; \alpha)$ and $\mathbf{H}^{+}(f, g; \alpha)$ be always homotopically equivalent (resp. isomorphic): when this is the case, we say that **A** *has homotopical homology* (resp. strict homotopical homology) and $\mathbf{H}(f, g; \alpha)$ denotes this common homotopy type (resp. isomorphism type). (Within IP4-homotopical categories, one can be more precise, and require that the canonical comparison $i$ be a homotopy equivalence, or an iso.) We shall see that the strict condition holds for chain complexes (3.4), and also in the unstable case of bounded chain complexes between fixed degrees (3.6).

To prove the existence of $i$, consider first that there is one map $x: Kg \to Kv_{\alpha}$ such that $kv_{\alpha}.x = cf.kg$, $\kappa v_{\alpha}.x = \kappa g$. Second, $xu_{\alpha} \simeq 0$, by the 2-dimensional universal property of the h-kernel $(Kv_{\alpha}, kv_{\alpha}, \kappa v_{\alpha})$; in fact, $kv_{\alpha}.xu_{\alpha} = cf.kg.u_{\alpha} = cf.f$, and the nullhomotopy $\gamma f: 0 \simeq cf.f$ is (strictly) coherent with $\kappa v_{\alpha}: 0 \simeq v_{\alpha}.kv_{\alpha}$, in the sense that $\kappa v_{\alpha}.xu_{\alpha} = \kappa g.u_{\alpha} = \alpha = v_{\alpha}.\gamma f$ (and $\kappa v_{\alpha}.0$ is the trivial endohomotopy). There is thus *some* homotopy $\rho: 0 \simeq xu_{\alpha}$ such that $kv_{\alpha}.\rho = \gamma f$. (In the IP4-homotopical case, we may determine $\rho$ via a *double homotopy* produced by the homotopy $\kappa v_{\alpha}.xu_{\alpha} = v_{\alpha}.\gamma f$ together with a connection of the cylinder or path functor.) Finally, this homotopy $\rho$ produces a map $i: Cu_{\alpha} \to Kv_{\alpha}$ such that $i.cu_{\alpha} = x$, $i.\gamma u_{\alpha} = \rho$, and (3) holds.

**2.6. A necessary condition.** If **A** has homotopical homology, then for every map $f$, $Cu_f$ and $Ku_f$ are contractible objects. In fact, applying 2.5.1 to the sequence $(f, cf; \gamma f)$, we get that $K(v_{\gamma f}) = K(1_{Cf}) = K(Cf)$ is a cocone; thus, the condition $Cu_f \simeq K(Cf)$ says that $Cu_f$ is contractible. In particular, if **A** has homotopical homology, the unit $u_X: X \to \Omega\Sigma X$ has a contractible h-kernel, while the counit $v_Y: \Sigma\Omega Y \to Y$ has a contractible h-cokernel.

This is sufficient to show that the homotopical category **Top•** of pointed topological spaces does *not* have homotopical homology: $Cu_X$ need not even be connected (taking, for instance, $X = \mathbf{S}^0$, the unit $\mathbf{S}^0 \to \Omega\mathbf{S}^1$ just reaches the trivial loop and the 'main' loop of the circle).

## 3. Homotopical exactness for chain complexes

This section deals with the usual homotopical structure of the category $Ch_*\mathbf{D}$ of unbounded chain complexes over an additive category **D**. Positive chain complexes are briefly considered (3.5-3.6); for bounded chain complexes, a notion of *weak exactness* is more adapted (3.7-9).

**3.1. Notation.** In the additive category **D**, a morphism $f: A \to B$ between two finite biproducts $A = \oplus A_i$, $B = \oplus B_j$ is determined by its components $f_{ji}: A_i \to B_j$ and can be written as a matrix $f = (f_{ji})$. But we prefer to write it as an 'expression in m variables'

(1) $f(a_1, a_2, ... a_m) = (\Sigma f_{1i} a_i, \Sigma f_{2i} a_i, ... \Sigma f_{ni} a_i)$,

as within a category of modules. (This can be formally justified letting $a_i$ denote the projection $A \to A_i$; $(a_1, a_2, ... a_m)$ is thus the identity of A, with the specification of the names of the projections.)

**3.2. The homotopical structure.** The category $Ch_*\mathbf{D}$ of (unbounded) chain complexes over the additive category **D** is equipped with the usual homotopies of chain maps

(1) $\alpha = (f, g, (\alpha_n)): f \simeq g: A \to B$,
$\alpha_n: A_n \to B_{n+1}$, $\qquad -f_n + g_n = \alpha_{n-1}\partial_n + \partial_{n+1}\alpha_n$.

These come from an IP4-homotopical structure, based on well known cylinder and path functors [G5, 6.5-6.8]. The opposite structure is isomorphic to $Ch_*(\mathbf{D}^{op})$, which is again of the previous type; this duality provides a choice of h-kernels and h-cokernels which reduces the structural isomorphisms to identities.

Let $f: A \to B$ be a map of chain complexes. The left-homotopical structure of $Ch_*\mathbf{D}$ is described as follows (as a mnemonic rule, note that a differential acting on a component shifted one degree up or down always requires a change of signs):

(2) $(Kf)_n = A_n \oplus B_{n+1}$, $\qquad \partial(a, b) = (\partial a, fa - \partial b)$,
$k: Kf \to A$, $k(a, b) = a$; $\qquad \kappa: 0 \simeq fk: Kf \to B$, $\kappa(a, b) = b$,

(3) $\Omega A = K(0 \to A)$, $\qquad (\Omega A)_n = A_{n+1}$, $\qquad \partial_n^{\Omega A} = -\partial_{n+1}$,
$\omega_A: 0 \simeq 0: \Omega A \to A$, $\qquad (\omega_A)_n = id(A_{n+1})$.

Analogously, the right-homotopical structure is given by:

(4) $(Cf)_n = A_{n-1} \oplus B_n$, $\qquad \partial(a, b) = (-\partial a, -fa + \partial b)$,
$c: B \to Cf$, $c(b) = (0, b)$; $\qquad \gamma: 0 \simeq cf: A \to Cf$, $\gamma(a) = (-a, 0)$,

(5) $\Sigma A \;=\; C(A \to 0)$, $\quad (\Sigma A)_n \;=\; A_{n-1}$, $\quad \partial_n^{\Sigma A} \;=\; -\,\partial_{n-1}$,

$\sigma_A\colon 0 \simeq 0\colon A \to \Sigma A$, $\quad (\sigma_A)_n \;=\; \mathrm{id}(A_{n-1})$.

**3.3. Stability.** This homotopical structure is strictly stable. To begin with, the functors $\Sigma, \Omega\colon \mathrm{Ch}_*\mathbf{D} \to \mathrm{Ch}_*\mathbf{D}$ (as realised above) are reciprocal automorphisms, and the unit and counit of the suspension-loop adjunction are identical.

Moreover, as to the stronger condition required in 2.2, an arbitrary map $f\colon A \to B$ produces a canonical homotopy $\rho\colon 0 \simeq 0\colon Kf \to Cf$ (2.2.2), which is represented by identities

(1) $\rho_n\colon (Kf)_n \to (Cf)_{n+1}$ $\qquad \rho_n(a, b) \;=\; (c\kappa - \gamma k)_n\,(a, b) \;=\; (0, b) + (a, 0) \;=\; (a, b)$,

so that the adjoint morphisms $V_f\colon \Sigma Kf \to Cf$ and $U_f\colon Kf \to \Omega Cf$ are identities as well.

Thus, $\mathrm{Ch}_*\mathbf{D}$ is also left and right *h*-stable (thm. 2.4); however, it is *not* left stable in the strict sense: the previous computations of h-kernels and h-cokernels show that $u_f\colon A \to Kcf$ is not an isomorphism, generally. By 2.4, the fibre-cofibre sequence of the chain morphism $f$ is strongly h-exact and sent to itself by the automorphisms $\Sigma$ and $\Omega$, with a three-place shift. Each h-exact sequence $(f, g; \alpha)$ is thus coherently equivalent (1.6) to a sequence $(B \to Cf \to \Sigma A; 0)$, which means a componentwise split short exact sequence of chain complexes (a notion which only requires the additivity of the basis $\mathbf{D}$); note that the homotopy $\alpha$ determines the splitting.

**3.4. Theorem: Homotopical homology of chain complexes.** In the category $\mathrm{Ch}_*\mathbf{D}$ of chain complexes over an additive category, the left and the right homotopical homology of an h-differential sequence (2.5) coincide, yielding one chain complex whose homotopy type vanishes iff the sequence is h-exact, or equivalently strongly h-exact

(1) $\mathbf{H}_n(f, g; \alpha) \;=\; X_{n-1}\oplus A_n\oplus Y_{n+1}$, $\qquad \partial(x, a, y) \;=\; (-\,\partial x, -\,fx + \partial a, -\,\alpha x + ga - \partial y)$.

**Proof.** Given the h-differential sequence $(f, g; \alpha)$, we show that $Cu_\alpha = Kv_\alpha$, a strict equality, using the construction of h-kernels and h-cokernels given in 3.2. The proof reduces to an easy computation; but it would be interesting to give a structural proof, in a more general situation.

(2)
$$\begin{array}{ccccc}
X & \xrightarrow{\;f\;} & A & \xrightarrow{\;g\;} & Y \\
{\scriptstyle u_\alpha}\downarrow & \nearrow{\scriptstyle k} & & \searrow{\scriptstyle c} & \uparrow{\scriptstyle v_\alpha} \\
Kg & \xrightarrow[cu_\alpha]{} & Cu_\alpha = Kv_\alpha & \xrightarrow[kv_\alpha]{} & Cf
\end{array}
\qquad (\alpha \text{ dashed from } X \text{ to } Y)$$

$\alpha\colon 0 \simeq gf\colon X \to Y$,

$k \;=\; kg$, $\quad \kappa \;=\; \kappa g\colon 0 \simeq g.k$,

$c \;=\; cf$, $\quad \gamma \;=\; \gamma f\colon 0 \simeq c.f$,

(3) $(Kg)_n \;=\; A_n\oplus Y_{n+1}$, $\qquad \partial(a, y) \;=\; (\partial a, ga - \partial y)$,

$k(a, y) \;=\; a$, $\quad \kappa(a, y) \;=\; y$, $\qquad u_\alpha(x) \;=\; (fx, \alpha x)$,

(4) $(Cf)_n \;=\; X_{n-1}\oplus A_n$, $\qquad \partial(x, a) \;=\; (-\,\partial x, -\,fx + \partial a)$,

$c(a) \;=\; (0, a)$, $\quad \gamma(x) \;=\; (-\,x, 0)$, $\qquad v_\alpha(x, a) \;=\; -\,\alpha x + ga$,

(5) $(Cu_\alpha)_n \;=\; (Kv_\alpha)_n \;=\; X_{n-1}\oplus A_n\oplus Y_{n+1}$,

$\partial(x, a, y) \;=\; (-\,\partial x, -\,fx + \partial a, -\,\alpha x + ga - \partial y)$,

$(cu_\alpha)(a, y) \;=\; (0, a, y)$, $\qquad (kv_\alpha)(x, a, y) \;=\; (x, a)$, $\qquad kv_\alpha.cu_\alpha \;=\; c.k$.

**3.5. Positive chain complexes.** Let the additive category $\mathbf{D}$ have *finite limits*. Then the category $Ch_p\mathbf{D}$ of *positive* chain complexes (null in negative degree) is homotopical and *left* h-stable. Indeed, it is a full subcategory of $Ch_*\mathbf{D}$, with reflector $L$ and coreflector $R$ ($L \dashv U \dashv R$; $LU = 1 = RU$)

(1) $(LA)_n = A_n \quad (n \geq 0),$

$(RA)_n = A_n \quad (n > 0), \qquad (RA)_0 = \mathrm{Ker}(\partial_0 \colon A_0 \to A_{-1}),$

and the (restricted) homotopy structure is produced by the composed adjunction

(2) $\mathbf{I} : Ch_p\mathbf{D} \rightleftarrows Ch_p\mathbf{D} : \mathbf{P}, \qquad \mathbf{I} = LIU \dashv \mathbf{P} = RPU.$

$Ch_p\mathbf{D}$ is closed in $Ch_*\mathbf{D}$ under h-cokernels (as $(Cf)_n = A_{n-1}\oplus B_n$ is null in negative degree), while the h-kernel of the chain morphism $f \colon A \to B$ in $Ch_p\mathbf{D}$ is given by the coreflector $R$ (applied to the h-kernel in $Ch_*\mathbf{D}$)

(3) $(Kf)_n = A_n\oplus B_{n+1} \quad (n > 0), \qquad (Kf)_0 =$ pullback of the pair $(f, \partial)$.

$Ch_p\mathbf{D}$ is left h-stable (2.3); this is proved directly by the computations below (3.6). But can also be deduced from the stability of $Ch_*\mathbf{D}$, applying a previous result on adjunctions (2.1) to the embedding $U^2 \colon (Ch_*\mathbf{D})^2 \to (Ch_p\mathbf{D})^2$, its adjoints and the left h-stable adjunction hcok/hker of $Ch_*\mathbf{D}$ (the hypothesis $U^2L^2cU^2 = cU^2$ follows from the closure of $Ch_p\mathbf{D}$ under h-cokernels).

On the other hand, $Ch_p\mathbf{D}$ is *not* right h-stable, as it follows immediately from its suspension and loop functors: they form a coreflexive adjunction $\Sigma \dashv \Omega$, $\Omega\Sigma = 1$

(4) $(\Sigma A)_n = A_{n-1},$

$(\Omega A)_n = A_{n+1} \quad (n > 0), \qquad (\Omega A)_0 = \mathrm{Ker}(\partial_1),$

(5) $\Omega\Sigma A = A,$

$(\Sigma\Omega A)_n = A_n \quad (n > 1), \qquad (\Sigma\Omega A)_1 = \mathrm{Ker}(\partial_1), \quad (\Sigma\Omega A)_0 = 0,$

which is not even a weak equivalence (3.7): $H_0(\Sigma\Omega A)$ is necessarily null.

Dually, if $\mathbf{D}$ has *finite colimits*, the category $Ch_n\mathbf{D}$ of negative chain complexes (null in positive degree) is homotopical and right h-stable. While, if $\mathbf{D}$ is finitely complete and cocomplete, the category $Ch_0^p\mathbf{D}$ of chain complexes concentrated in degrees $n \in [0, p]$ $(p \geq 1)$ is *h-semistable*.

**3.6. Computations.** For positive chain complexes, the left h-stable property can be easily checked directly: given $f \colon X \to A$, we show that $u_f \colon X \to K$ is a *deformation retract*: (where $K = Kcf$ and $Cf$ is computed as in 3.4.4)

$$(1) \qquad \begin{array}{ccccc} X & \xrightarrow{f} & A & \xrightarrow{cf} & Cf \\ {\scriptstyle u}\downarrow\uparrow{\scriptstyle u'} & \nearrow_{k} & & & \\ K & & & & \end{array} \qquad \begin{array}{l} \gamma = \gamma f \colon 0 \simeq cf.f \colon X \to Cf \\ \kappa = \kappa cf \colon 0 \simeq cf.k \colon K \to Cf \\ u = u_f, \end{array}$$

(2) $K_n = A_n\oplus X_n\oplus A_{n+1} \quad (n > 0), \qquad K_0 =$ pullback of $(cf_0, \partial_1) \subset A_0\oplus X_0\oplus A_1,$

$\partial(a, x, a') = (\partial a, \partial x, a - \partial a' + fx),$

$k(a, x, a') = a, \quad \kappa(a, x, a') = (x, a'), \qquad u(x) = (fx, -x, 0),$

letting $u'(a, x, a') = -x$. It is a chain map, with $u'u = 1$ and there is a homotopy $\sigma\colon uu' \simeq 1_K$

(3) $\sigma_n\colon K_n \to K_{n+1}$, $\qquad\qquad \sigma(a, x, a') = (a', 0, 0)$.

As to the complex of homotopical homology (3.4), it is interesting to note that the relation $Cu_\alpha = Kv_\alpha$ providing $\mathbf{H}(f, g; \alpha)$ holds also in the present case of *positive* chain complexes, with

(4) $\mathbf{H}_0(f, g; \alpha) = (Kg)_0 = (Cu_\alpha)_0 = (Kv_\alpha)_0 = \mathrm{pb}(g_0\colon A_0 \to Y_0 \leftarrow Y_1 :\partial_1)$.

We already know (by 2.5) that *the h-differential sequence* $(f, g; \alpha)$ *is h-exact iff it is left h-exact, iff* $\mathbf{H}(f, g; \alpha)$ *is contractible*.

But being *right* h-exact is indeed a stronger condition (even up to *weak equivalence*). For instance, the typical left h-exact sequence $(\Omega A \to 0 \to A; \omega_A)$ is right h-exact iff the counit $v_A\colon \Sigma\Omega A \to A$ is a homotopy equivalence, which requires that $H_0(A) = H_0(\Sigma\Omega A) = 0$. The homotopical homology $\mathbf{H}(\Omega A, A; \omega_A)$ of that sequence consists of the h-cokernel of $1\colon \Omega A \to \Omega A$, i.e. the cone $C(\Omega A)$, which is always contractible (independently of right h-exactness, of course).

Dually, for negative chain complexes, $\mathbf{H}(f, g; \alpha)$ is contractible iff the sequence is h-exact, iff it is right h-exact (while left h-exactness is a stronger condition); its 0-component is computed as a pushout

(5) $\mathbf{H}_0(f, g; \alpha) = (Cf)_0 = (Cu_\alpha)_0 = (Kv_\alpha)_0 = \mathrm{po}(X_{-1} \leftarrow X_0 \to A_0)$;

Finally, in $\mathrm{Ch}_0^p\mathbf{D}$ (for $\mathbf{D}$ finitely complete and cocomplete), the object $\mathbf{H}(f, g; \alpha) = Cu_\alpha = Kv_\alpha$ is still well-defined, and computed by a limit or colimit in the boundary degrees

(6) $\mathbf{H}_0(f, g; \alpha) = \mathrm{pb}(A_0 \to Y_0 \leftarrow Y_1)$, $\qquad\qquad \mathbf{H}_p(f, g; \alpha) = \mathrm{po}(X_{p-1} \leftarrow X_p \to A_p)$.

For the particular case $p = 1$, the homotopy type of $\mathbf{H}(f, g; \alpha)$ (a complex reduced to a map!) still measures the default of h-exactness of the sequence, as we prove in the next section (4.5). For $p > 1$, this holds in a weaker sense, linked to *weak equivalences*, the chain maps which induce an isomorphism in homology, as we prove now.

**3.7. Weak equivalences.** Let $\mathbf{D}$ be abelian. The homology functors $H_n\colon \mathbf{A} = \mathrm{Ch}_0^p\mathbf{D} \to \mathbf{D}$ $(0 \leq n \leq p)$ satisfy the following 'self-dual axioms', combining the ones of homology and homotopy theories

(i) *homotopy invariance*: if $f \simeq g$ in $\mathbf{A}$, then $H_n(f) = H_n(g)$;

(ii) *exactness*: for every map $f\colon A \to B$ in $\mathbf{A}$, the sequences

(1) $H_0(Kf) \to H_0(A) \to H_0(B)$, $\qquad\qquad H_p(A) \to H_p(B) \to H_p(Cf)$,

are exact in $\mathbf{D}$;

(iii) *stability*: there are two natural isomorphisms, $\omega\colon H_n(\Omega A) \to H_{n+1}(A)$ $(0 \leq n < p)$ and $\sigma\colon H_n(A) \to H_{n+1}(\Sigma A)$ $(0 \leq n < p)$, which are coherent with the morphisms $H_n \to H_n\Omega\Sigma$ $(n < p)$ and $H_m\Sigma\Omega \to H_m$ $(m > 0)$ induced by the unit and counit of the adjunction $\Sigma \dashv \Omega$, forming commutative diagrams (of isomorphisms of abelian groups), for $n < p$, $m > 0$:

(2)
$$\begin{array}{ccc} H_n(A) & \longrightarrow & H_n(\Omega\Sigma A) \\ & {\scriptstyle \sigma A}\searrow \quad \swarrow{\scriptstyle \omega\Sigma A} & \\ & H_{n+1}(\Sigma A) & \end{array} \qquad\qquad \begin{array}{ccc} H_m(\Sigma\Omega A) & \longrightarrow & H_m(A) \\ & {\scriptstyle \sigma\Omega A}\nwarrow \quad \nearrow{\scriptstyle \omega A} & \\ & H_{m-1}(\Sigma A) & \end{array}$$

As a consequence, every map f has two exact **D**-sequences, the H-*fibre* and H-*cofibre* sequence

(3) $0 \to H_p(Kf) \to H_p(A) \to H_p(B) \to \ ... \ \to H_0(Kf) \to H_0(A) \to H_0(B),$

(4) $H_p(A) \to H_p(B) \to H_p(Cf) \to \ ... \ \to H_0(A) \to H_0(B) \to H_0(Cf) \to 0.$

For the first sequence, it suffices to apply $H_0$ to the fibre sequence of f, taking into account that $H_n(A) \cong H_0(\Omega^n A)$ for $0 \leq n \leq p$ and that $\Omega^{p+1} = 0$; exactness follows from the properties of the fibre sequence (1.8), together with the axioms of exactness and homotopy invariance for $H_0$.

Now, say that a map $f: A \to B$ is a *weak equivalence* (with respect to the theory $H_*$) if $H_n(f)$ is iso for all n; similarly, the object A is *weakly null* (or weakly contractible) if $H_n(A) = 0$ for all n. Plainly, by the usual properties of exact sequences in abelian categories:

(a) Kf is weakly null if and only if $H_n(f)$ is iso for $n > 0$ and monic for $n = 0$,

(a*) Cf is weakly null if and only if $H_n(f)$ is iso for $n < p$ and epi for $n = p$.

Thus, for $p = 2$, the complex $A = (A_2 \to A_1 \to A_0)$ is a weakly null object iff it forms a short exact sequence; but it is easy to see that it is nullhomotopic iff it forms a *split*short exact sequence.

More generally, all this still hold for any pointed homotopical category **A** with $\Omega^{p+1} = 0 = \Sigma^{p+1}$, equipped with a theory $H_n$: **A** $\to$ **D** with values in an abelian category and satisfying the axioms above for $0 \leq n \leq p$. This setting can be easily adapted to the positive or unbounded cases.

**3.8. Lemma.** In the general setting just mentioned, let an h-differential sequence $(f, g; \alpha)$ be given. Then the map $i: Cu_\alpha \to Kv_\alpha$ (2.5) is a weak equivalence.

**Proof.** The commutative diagram 2.5.1 gives a commutative diagram, in the abelian category **D**

(1)

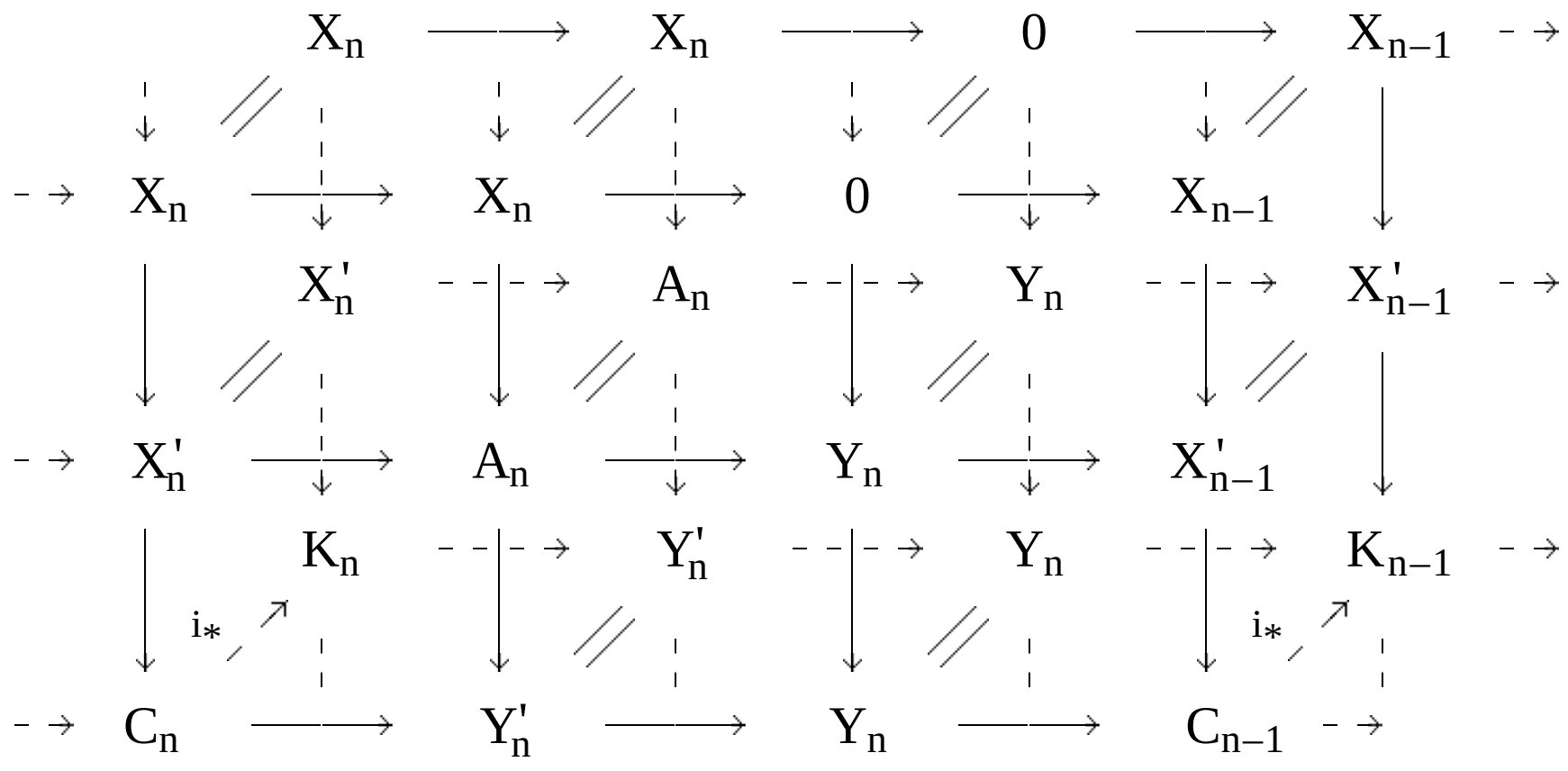

with $C_n = H_n(Cu_\alpha)$ and $K_n = H_n(Kv_\alpha)$. Its rows and columns are exact, except - possibly - the lowest front row and the left-hand back column, which are just known to be of order two. Moreover, in each row or column, each term $C_n$ or $K_n$ has at least two terms on each side (possibly null). By a sort of variation of the 3×3-Lemma, it follows that $i_{*n} = H_n(i)$ is iso. (Of course, one only needs to verify this for modules, by diagram chasing.)

**3.9. Theorem and Definition: Weak exactness.** In the general setting mentioned at the end of 3.7, let an h-differential sequence $(f, g; \alpha)$ be given, and consider again the diagram 1.3.1. The following conditions are equivalent:

(a) the map $v: Cf \to Ckg$ is a weak equivalence (with respect to $H_*$),

(b) the morphism $H_n(u_\alpha): H_n(X) \to H_n(Kg)$ is iso for $n < p$, and epi for $n = p$,

(c) the object $\mathbf{H}^-(f, g; \alpha) = Cu_\alpha$ is weakly null,

(a') the map $u: Kcf \to Kg$ is a weak equivalence (with respect to $H_*$),

(b') the morphism $H_n(v_\alpha): H_n(Cf) \to H_n(Y)$ is iso for $n > 0$, and monic for $n = 0$,

(c') the object $\mathbf{H}^+(f, g; \alpha) = Kv_\alpha$ is weakly null.

If they are satisfied, the sequence $(f, g; \alpha)$ will be said to be *w-exact*. This notion is strictly weaker than h-exactness, provided that $p \geq 2$ and $\mathbf{D}$ has some non-split short exact sequence.

**Proof**. The preceding lemma shows that (c) and (c') are equivalent, while the equivalence of (b) and (c) has been considered at the end of 3.7. Therefore, we only need to prove that (a) and (b) are equivalent. This follows from the Five Lemma in $\mathbf{D}$: take the diagram

$$(1) \qquad \begin{array}{ccccc} Kcf & \xrightarrow{kcf} & A & \xrightarrow{cf} & Cf \\ {\scriptstyle u}\downarrow & & \| & & \downarrow{\scriptstyle v_\alpha} \\ Kg & \xrightarrow{kg} & A & \xrightarrow{g} & Y \end{array}$$

and apply the (natural) H-fibre sequence, together with the fact that $H_0(cf)$ is epi (3.7.4)

$$(2) \qquad \begin{array}{ccccccccccc} 0 & \to & H_p(Kcf) & \longrightarrow & H_p(A) & \ldots & \longrightarrow & H_0(A) & \longrightarrow & H_0(Cf) & \longrightarrow & 0 \\ & & {\scriptstyle u_*}\downarrow & & \| & & & \| & & {\scriptstyle v_*}\downarrow & & \downarrow \\ 0 & \to & H_p(Kg) & \longrightarrow & H_p(A) & \ldots & \longrightarrow & H_0(A) & \longrightarrow & H_0(Y) & \longrightarrow & Cok(H_0 g) \end{array}$$

Finally, it suffices to consider a sequence $0 \to A \to 0$, which is h-exact (resp. w-exact) iff $A$ is nullhomotopic (resp. w-null), which - under our assumptions - are indeed non-equivalent conditions by a remark at the end of 3.7.

## 4. Homotopical categories of morphisms

The case $Ch_0^1\mathbf{D} = \mathbf{D}^2$ gives a simple *h-semistable* category having homotopical homology.

**4.1. The homotopy structure.** In this section, $\mathbf{D}$ is always an additive category with finite limits and colimits and $\mathbf{A} = Ch_0^1\mathbf{D} = \mathbf{D}^2$ is the category of arrows of $\mathbf{D}$, enriched with the two dimensional structure of chain complexes (in degrees 0, 1), already considered in the previous section (3.5-6): a pointed homotopical, *h-semistable* category.

An object will be written as $A = \partial_A: A' \to A''$ and a map as $f = (f', f''): A \to B$. A nullhomotopy $\alpha: 0 \simeq f$ is determined by a *diagonal* map $\alpha: A'' \to B'$

$$(1) \qquad \begin{array}{ccc} A' & \xrightarrow{\ f'\ } & B' \\ \partial\downarrow & \overset{\alpha}{\nearrow} & \downarrow\partial \\ A'' & \xrightarrow[\ f''\ ]{} & B'' \end{array} \qquad\qquad f' \ = \ \alpha\partial, \qquad f'' \ = \ \partial\alpha.$$

while a homotopy $\alpha: f \simeq g: A \to B$ is given by a nullhomotopy $\alpha: 0 \simeq -f + g$ $(-f' + g' = \alpha\partial,\ -f'' + g'' = \partial\alpha)$. As an exception in the family $\mathrm{Ch}_0^p\mathbf{D}$, *we have here a* 2*-category* (essentially, because two-dimensional homotopies vanish, being produced by morphisms of degree 2): given two consecutive nullhomotopies $\alpha: 0 \simeq f$ and $\beta: 0 \simeq g$, the reduced interchange holds, $g\alpha = \beta f$ (both homotopies are determined by the same diagonal, $g'\alpha = \beta\partial\alpha = \beta f'': A'' \to C'$).

Plainly, an object $\partial_A: A' \to A''$ is contractible iff $\partial_A$ is an isomorphism. As a consequence (for a non trivial $\mathbf{D}$), a homotopy equivalence need not be an iso; however, there is a relevant case where the two notions coincide.

**4.2. A coherence lemma.** If, in the morphism $f: A \to B$, either $f'$ or $f''$ is an iso, then $f$ is a homotopy equivalence iff it is an isomorphism.

**Proof.** Let $f' = 1_Z$ $(Z = A' = B')$ and choose an adjoint equivalence $A \rightleftarrows B$

$$(1) \qquad \begin{array}{ccccccc} Z & = & Z & \xrightarrow{\ g\ } & Z & = & Z \\ \partial\downarrow & \overset{\alpha}{\swarrow} & \downarrow\partial & \nearrow & \downarrow\partial & \overset{\beta}{\nearrow} & \downarrow\partial \\ A'' & \xrightarrow[f'']{} & B'' & \xrightarrow[g'']{} & A'' & \xrightarrow[f'']{} & B'' \end{array} \qquad \begin{array}{l} \alpha: 1_A \simeq gf, \quad \beta: fg \simeq 1_B \\ \\ f\circ\alpha + \beta\circ f \ = \ 0, \quad \alpha\circ g + g\circ\beta \ = \ 0. \end{array}$$

Then $f''$ admits an inverse mapping: $h = g'' + \partial_A\beta$. In fact, $f''h = f''g'' + f''\partial_A\beta = 1 - \partial_B\beta + \partial_B\beta = 1$ and $hf'' = g''f'' + \partial_A\beta f'' = 1 + \partial_A\alpha - \partial_A\alpha = 1$ (use the triangle identity $f\circ\alpha + \beta\circ f = 0$). Similarly, if $f'' = 1_Z$ $(Z = A'' = B'')$, then $f'$ has an inverse mapping, namely $g' - \alpha\partial_B: B' \to A'$.

**4.3. h-kernels and h-cokernels.** If $f: A \to B$ is the commutative square 4.1.1, form the pullback $K$ of $(f'', \partial_B)$ and the pushout $C$ of $(\partial_A, f')$; this yields an internal bicartesian square

$$(1) \qquad \begin{array}{ccccccccc} & \overset{\partial}{\nearrow} & & \overset{k''}{\nearrow} & A'' & \overset{\gamma}{\searrow} & & \overset{f''}{\searrow} & \\ A' & \xrightarrow{\ k'\ } & K & & & & C & \xrightarrow{\ c''\ } & B'' \\ & \underset{f'}{\searrow} & & \underset{\kappa}{\searrow} & B' & \underset{c'}{\nearrow} & & \underset{\partial}{\nearrow} & \end{array}$$

from which the h-kernel and the h-cokernel of $f$ are constructed

$$(2) \qquad \begin{array}{ccccccc} A' & = & A' & \xrightarrow{\ f'\ } & B' & \xrightarrow{\ c'\ } & C \\ k'\downarrow & \swarrow & \downarrow & \overset{\kappa}{\nearrow} & \downarrow & \overset{\gamma}{\nearrow} & \downarrow c'' \\ K & \xrightarrow[k'']{} & A'' & \xrightarrow[f'']{} & B'' & = & B'' \end{array} \qquad \begin{array}{l} Kf \ = \ k', \quad kf \ = \ (1, k''), \quad \kappa f \ = \ \kappa, \\ \\ Cf \ = \ c'', \quad cf \ = \ (c', 1), \quad \gamma f \ = \ \gamma. \end{array}$$

**4.4. Left and right exactness.** Given an h-differential sequence (f, g; α), the properties of h-exactness are characterised as follows, up to coherent equivalence (1.6)

(1) 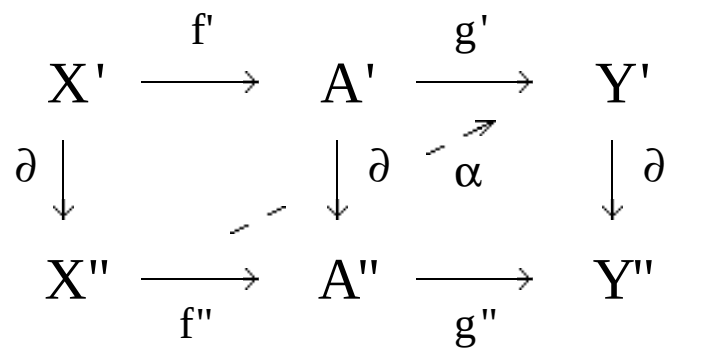

$\alpha: 0 \simeq gf: X \to Y$

(a) in the typical left h-exact sequence, (X'', α, f'') is the pullback of $(g, \partial_Y)$ and $f' = 1_{A'}$,

(a*) in the typical right h-exact sequence, (Y', α, g') is the pushout of $(f', \partial_X)$ and $g'' = 1_{A''}$,

(b) the typical strongly h-exact sequence is of the following type, determined by an arbitrary factorisation $\partial_A = h''h'$

(2) 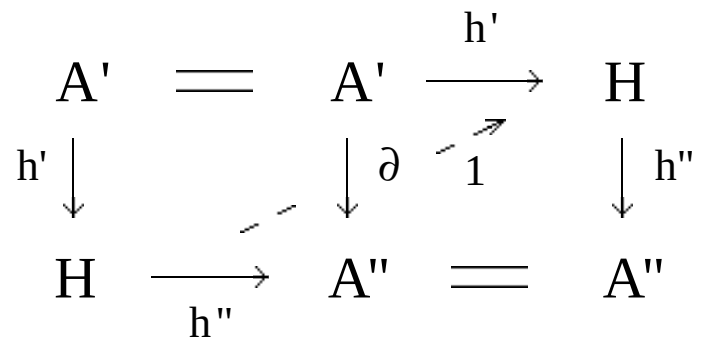

**4.5. Homotopical homology.** We already know, from 3.6, that an h-differential sequence α: 0 ≃ gf always has homotopical homology; we prove now that the homotopy type of the latter does measure the default of h-exactness.

First, to compute **H**(f, g; α), let us form in **D** the left-hand commutative diagram below, where (B, c, γ) is the pushout of $(\partial_X, f')$, while (Z, k, κ) is the pullback of $(\partial_Y, g'')$ and w: B → Z is the induced morphism

(1) 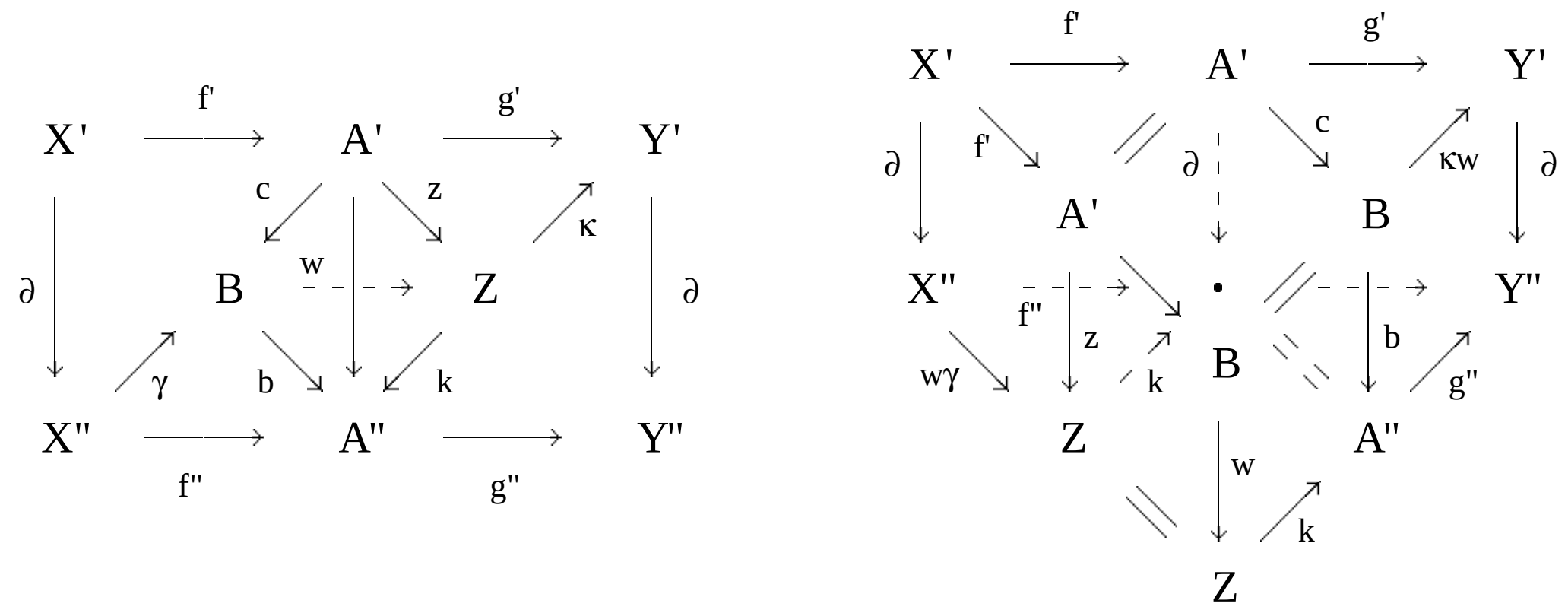

Now, the right-hand diagram shows that the object w: B → Z is $\mathbf{H}(f, g; \alpha) = Cu_\alpha = Kv_\alpha$

(2) $\mathrm{hker}(g) = (z, (1, k): z \to \partial_A, \kappa)$, $\qquad \mathrm{hcok}(f) = (b, (c, 1): \partial_A \to b, \gamma)$,

$u_\alpha = (f', w\gamma): X \to z$, $\qquad v_\alpha = (\kappa w, g''): b \to Y$,

$\mathrm{hcok}(u_\alpha) = (w, (c, 1): z \to w, \gamma)$, $\qquad \mathrm{hker}(v_\alpha) = (w, (1, k): w \to b, \kappa)$.

Finally, the diagram 1.3.1 becomes

(3)

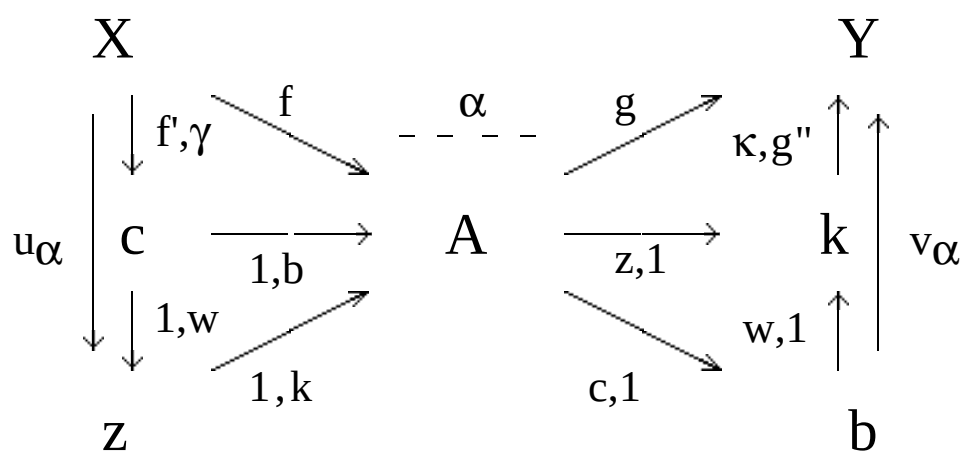

$u_\alpha = (f', w\gamma)$

$v_\alpha = (\kappa w, g")$

therefore, our sequence is h-exact iff $k(v_\alpha) = (1_{A'}, w)$ is a homotopy equivalence, iff (by lemma 4.2) $w$ is an isomorphism of **D**, iff **H**(f, g; α) is contractible (4.1), iff **H**(f, g; α) is weakly null (3.7), iff the sequence is w-exact (3.9).

## References

[Be] R. Betti, *Adjointness in descent theory*, J. Pure Appl. Algebra **116** (1997), 41-47.

[Do] A. Dold, *Partitions of unity in the theory of fibrations*, Ann. of Math. **78** (1963), 223-255.

[F1] P. Freyd, *Representations in abelian categories*, in: Proc. Conf. Categ. Algebra, La Jolla, 1965, Springer, Berlin 1966, pp. 95-120.

[F2] P. Freyd, *Stable homotopy*, in: Proc. Conf. Categ. Algebra, La Jolla, 1965, Springer, Berlin 1966, pp. 121-176.

[F3] P. Freyd, *Stable homotopy II*, in: Proc. Symp. Pure Maths. XVII, Amer. Math. Soc., Providence, RI, 1970, pp. 161-183.

[G1] M. Grandis, *On the categorical foundations of homological and homotopical algebra*, Cahiers Top. Géom. Diff. Catég. **33** (1992), 135-175.

[G2] M. Grandis, *Homotopical algebra in homotopical categories*, Appl. Categ. Structures **2** (1994), 351-406.

[G3] M. Grandis, *Cubical monads and their symmetries*, in: Proc. of the Eleventh Intern. Conf. on Topology, Trieste 1993, Rend. Ist. Mat. Univ. Trieste, **25** (1993), 223-262.

[G4] M. Grandis, *Homotopical algebra and triangulated categories,* Math. Proc. Cambridge Philos. Soc. **118** (1995), 259-295.

[G5] M. Grandis, *Categorically algebraic foundations for homotopical algebra*, Appl. Categ. Structures **5** (1997), 363-413.

[G6] M. Grandis, *Weak subobjects and the epi-monic completion of a category*, J. Pure Appl. Algebra, to appear.

[Ha] R. Hartshorne, *Residues and duality*, Lecture Notes in Math. vol. 20, Springer, Berlin 1966.

[HK1] K. Hardie - K.H. Kamps, *Track homotopy over a fixed space*, Glasnik Matematički **24** (1989), 161-179.

[HK2] K. Hardie - K.H. Kamps, *Coherent homotopy over a fixed space*, in: Handbook of Algebraic Topology, Elsevier, Amsterdam, 1995, 195-211.

[HKP] K. Hardie - K.H. Kamps - T. Porter, *The coherent homotopy category over a fixed space is a category of fractions*, Topology and its Appl. **40** (1991), 265-274.

[Ka] D. M. Kan, *Functors involving c.s.s. complexes*, Trans. Amer. Math. Soc. 87 (1958), 330-346.

[KV] S. Kasangian - E.M. Vitale, *Factorization systems for symmetric cat-groups*, Th. Appl. Categories **7**-5 (2000), 47-70.

[P1] D. Puppe, *Homotopiemengen und ihre induzierten Abbildungen I*, Math. Z. **69** (1958), 299-344.

[P2] D. Puppe, *On the formal structure of stable homotopy theory*, in: Colloquium on Algebraic Topology. Mat. Inst., Aarhus Univ., 1962, 65-71.

[P3] D. Puppe, *Stabile Homotopietheorie I*, Math. Ann. **169** (1967), 243-274.

[Ve] J. L. Verdier, *Catégories dérivées*, in: Séminaire de Géométrie algébrique du Bois Marie SGA 4 $^1/_2$, Cohomologie étale, Lecture Notes in Math. vol. 569, Springer, Berlin 1977, pp. 262-311.

[Vo] R.M. Vogt, *A note on homotopy equivalences*, Proc. Amer. Math. Soc. **32** (1972), 627-629.